\documentclass[11pt]{article}
\usepackage{amsmath,amsthm}

\def\ftoday{le \space\number\day \space\ifcase\month\or
  janvier\or f\'evrier\or mars\or avril\or mai\or juin\or
  juillet\or ao\^ut\or septembre\or octobre\or novembre\or d\'ecembre\fi
  \space\number\year}

\def\real{I\kern-0.20em R}

\def\integer{I\kern-0.20em N}
\def\relative{{\rm \rlap Z\kern 2.2pt Z}}
\def\cc{\kern-.25em{\c c}}

\def\bc{\begin{center}}
\def\ec{\end{center}}
\def\=def{\stackrel{{\rm def}}{=}}

\newcommand\vfo[2]{{\cal X}_{#1}^{#2}}
\newcommand\fvfo[2]{\widehat{\cal X}_{#1}^{#2}}

\newcommand\pvf[3]{{\cal P}_{#1}^{#2,#3}}
\newcommand\hvf[2]{{\cal P}_{#1}^{#2}}
\newcommand\lie[1]{{\frak #1}}
\renewcommand\hom[1]{\text{Hom}_{\Bbb #1}}

\newcounter{indconst}
\newcounter{auxconst}

\def\bit{\begin{itemize}}
\def\eit{\end{itemize}}
\def\ben{\begin{enumerate}}
\def\een{\end{enumerate}}
\def\bde{\begin{description}}
\def\ede{\end{description}}
 
\def\beq{\begin{equation}}
\def\eeq{\end{equation}}
\def\bfi{\begin{figure}[hbt] \begin{center}}
\def\efi{\end{center} \end{figure}}

\def\bce{\begin{center}}
\def\ece{\end{center}}

\newcommand{\Preuve}{\medskip\noindent{\underline{Preuve:}}\quad\medskip}        
\newcommand{\cqfd  }{\hfill$\square$}
 
\newtheorem {theo} {Th\'eor\`eme}[section]

\newtheorem {coro} {Corollaire} [section]

\newtheorem {lemm} {Lemme}[section]
\newtheorem {lemms} {Lemma}[section]
\newtheorem {prop} {Proposition}[section]

\newtheorem {defi} {D\'efinition}[section]

\newtheorem {rem}{Remarque} [section]

\newtheorem {ex} {Exemple}[section]

\input amssym.def
\input amssym

\def\abstractname{R\'{e}sum\'{e}}

\begin{document}

\title{Normalisation holomorphe d'alg\`ebres de type Cartan de champs de vecteurs holomorphes singuliers}
\author{Laurent Stolovitch \thanks{CNRS UMR 5580, Laboratoire Emile Picard,
Universite Paul Sabatier, 118 route de Narbonne,
31062 Toulouse cedex 4, France. Courriel : {\tt stolo@picard.ups-tlse.fr}} }
\date{\ftoday}
\maketitle
\begin{abstract}
Nous montrons r\'esultat de normalisation holomorphe d'une famille commutative de champs de vecteurs holomorphes 
au voisinage de leur point singulier. Pour ce faire, nous supposons qu'une condition diophantienne portant sur
une alg\`ebre de Lie commutative de champs lin\'eaires associ\'ee aux parties lin\'eaires 
de la famille est satisfaite. D'autre part, nous imposons certaines conditions alg\'ebriques sur leur forme normale formelle. 
Les champs d'une telle famille, sauf un, peuvent ne pas avoir de partie lin\'eaire \`a l'origine.
\end{abstract}
\def\abstractname{Abstract}
\begin{abstract}
We consider a commutative familly of holomorphic vector fields in an neighbourhood of 
a common singular point, say $0\in \Bbb C^n$. Let $\lie g$ be a commutative complex Lie algebra 
of dimension $l$. Let $\lambda_1,\ldots,\lambda_n\in \lie g^*$ and 
let us set $S(g)=\sum_{i=1}^n\lambda_i(g)x_i\frac{\partial}{\partial x_i}$. We assume that this Lie morphism 
is {\bf diophantine} in the sense that a diophantine condition $(\omega(S))$ 
is satisfied. Let $X_1$ be a holomorphic vector fields in a neighbourhood of $0\in \Bbb C^n$. We assume that 
its linear part $s$ is regular relatively to $S$, that is belongs to $S(\lie g)$ and has the same formal commutator 
as $S$. Let $X_2,\ldots, X_l$ be holomorphic vector fields vanishing at $0$ and commuting with $X_1$.
Then there exists a formal diffeomorphism of $(\Bbb C^n,0)$ such that the familly of vector fields are in {\bf normal form} in 
these formal coordinates. This means that each element of the familly commutes with $s$. 
We show that, if the normal forms of the $X_i$'s belongs to $\widehat {\cal O}_n^S\otimes S(\lie g)$ 
($\widehat {\cal O}_n^S$ is the ring of formal first integral of $S$) and their junior parts are free over $\widehat {\cal O}_n^S$, then 
there exists a holomorphic diffeomorphism of $(\Bbb C^n,0)$ which transforms the familly into a normal form.
The elements of the familly, but one, may not have a non-zero linear part.
\end{abstract}
\tableofcontents
\section{Introduction}

Dans cet article, nous nous int\'eressons au comportement local des trajectoires de champs de vecteurs 
holomorphes au voisinage d'un point singulier commun et plus pr\'ecis\'ement, nous nous 
int\'eressons au probl\`eme de la normalisation holomorphe de ces champs. Depuis Poincar\'e, de nombreux 
travaux y ont \'et\'e consacr\'e; en particulier, ceux de V.I. Arnold, A.D. Bruno, J. Ecalle, J.-P. Ramis, l'auteur.

\subsection{Rappels historiques}

Commen\c{c}ons par une id\'ee simple. Supposons que l'on veuille \'etudier l'orbite $\{A^kx\}_{k\in \Bbb N}$ d'un point $x\in \Bbb C^n$ par l'action 
d'une matrice $A\in {\cal M}_n(\Bbb C)$. Un des moyens \'efficaces de proc\'eder est de transformer la matrice $A$ en une {\bf forme de Jordan}.
A partir de cette matrice de Jordan et de la matrice de passage, toute la dynamique peut \^etre comprise. 
L'id\'ee de Poincar\'e est de proc\'eder de m\^eme pour les champs de vecteurs singuliers. Le but est alors de transformer un champs de vecteurs s'annulant 
en $0\in \Bbb C^n$ en un champs de vecteurs "plus simple" par un diff\'eomorphisme local de $(\Bbb C^n,0)$ fixant l'origine  et
tangent \`a l'identit\'e en ce point. On esp\`ere alors pouvoir \'etudier la dynamique du mod\`ele simple puis 
"remonter" ces informations au syst\`eme de d\'epart via la transformation.
Si $X$ est un tel champ de vecteurs et $\phi$ un diff\'eomorphisme de $(\Bbb C^n,0)$, 
alors l'action de $\phi$ sur $X$ est le champ de vecteurs $X'$ d\'efinit par $X'(\phi(x))=D\phi(x)X(x)$. 
Supposons pour simplifier que la partie lin\'eaire de $X$ soit un champ de vecteurs diagonal, c'est-\`a-dire 
$DX(0)(x)=\sum_{i=1}^n\lambda_i x_i\frac{\partial}{\partial x_i}\=def S$ o\`u les $\lambda_i$ sont des nombres complexes non tous nuls.
On montre alors \cite{Arn2,roussarie-ast,Camacho-port} qu'il existe un diff\'eomorphisme formel $\hat\Phi$ (qui n'est pas unique) tel que 
$$
\hat\Phi^*X=\sum_{i=1}^n{\lambda_ix_i\frac{\partial}{\partial x_i}}+\sum_{i=1}^n{\left(\sum_{(Q,\lambda)=\lambda_i} a_{i,Q} x^Q \right)\frac{\partial}{\partial x_i}}
$$
o\`u la somme porte sur les muliindices $Q\in \Bbb N^n$, $|Q|\geq 2$ et l'indice $i$ qui satisfont 
$(Q,\lambda)\=def \sum_{j=1}^n q_j\lambda_j=\lambda_i$ 
et o\`u les $a_{i,Q}$'s sont de nombres complexes. Comme d'habitude, si $Q=(q_1,\ldots,q_n)$, alors $x^Q=x_1^{q_1}\ldots x_n^{q_n}$ 
et $|Q|=q_1+\cdots q_n$. 
Le champ de vecteur (formel) auquel on arrive est appel\'e {\bf forme normale} formelle (de Poincar\'e-Dulac). 
Elle n'est en g\'en\'eral pas unique. On peut dire que c'est une "forme de Jordan" de $X$ dans le sens o\`u elle s'\'ecrit 
$\hat\Phi^*X=S + N$, $N$ \'etant le champ de vecteurs "nilpotent" d\'efini par les s\'eries et v\'erifiant $[S,N]=0$. Ici, 
$[.,.]$ est le crochet de Lie des champs de vecteurs. On voit appara\^{\i}tre les obstructions formelles \`a la lin\'earisation :
ceux sont les $Q\in \Bbb N^n$ et $1\leq i\leq n$ qui v\'erifient la {\bf relation de r\'esonance} $(Q,\lambda)=\lambda_i$. Ainsi, lorsque les valeurs propres 
$\lambda_1,\ldots, \lambda_n$ ne sont pas r\'esonantes, il existe un diff\'eomorphisme formel (unique dans ce cas) 
qui lin\'earise le champ de vecteurs $X$. Bien $X$ soit holomorphe au voisinage de $0\in \Bbb C^n$, il n'en demeure pas moins que le diff\'eomorphisme 
peut ne pas \^etre holomorphe au voisinage de l'origine. Un condition suffisante pour assurer l'holomorphie de ce diff\'eomorphisme 
est la condition diophantienne $(\omega)$ de Bruno portant sur les valeurs propres de $S$. 
Nous exprimerons cette condition plus loin; n\'eanmoins, pour donner une id\'ee au lecteur, cette condition 
est plus faible que la condition de C.L. Siegel \cite{Siegel,Arn2} : il existe $C>0$ et $\mu\geq 0$ tels que 
pour tout $Q\in \Bbb N^n$, $|Q|\geq 2$, et $1\leq i\leq n$, $\left|(Q,\lambda)-\lambda_i\right|\geq \frac{C}{|Q|^{\mu}}$.
Sous une telle hypoth\`ese, le champ de vecteurs $X$ est holomorphiquement lin\'earisable. En dimension $2$, 
J.-C. Yoccoz et R. P\'erez-Marco ont d\'emontr\'e la n\'ec\'essit\'e de la condition de Bruno \cite{perez-yoccoz} 
\`a partir de la n\'ec\'essit\'e de la condition de Bruno pour lin\'eariser un diff\'eomorphisme de $(\Bbb C,0)$ \cite{Yoccoz,Yoccoz-ast}.

Comme nous l'avons vu, l'holomorphie d'une transformation normalisante est tr\`es li\'ee au probl\`eme des 
{\bf petits diviseurs}. Elle est aussi li\'ee, dans le cas r\'esonnant, \`a l'existence d'int\'egrale premi\`ere formelle non-triviale.
Par exemple, une forme normale du champ $x^2\frac{\partial}{\partial x}+(x+y)\frac{\partial}{\partial y}$ est 
$x_1^2\frac{\partial}{\partial x_1}+y_1\frac{\partial}{\partial y_1}$. Le diff\'eomorphisme normalisant est, dans ce cas, 
$x=x_1,\;y=y_1+\sum_{k\geq 1}(k-1)!x_1^k$. Il n'est clairement pas holomorphe au voisinage de $0\in \Bbb C^2$.
Cette situation a engendr\'e de nombreux travaux \cite{Malg-diffeo, Ram-Mart1, Ram-Mart2, Voronin, Ecalle-mart, Stolo-clas}.

Dans quelle situation, un champ holomorphe singulier est-il holomorphiquement normalisable? 
A.D. Bruno a donn\'e une reponse a cette question dans un imposant m\'emoire \cite{Bruno}; un de ces r\'esultat est le suivant~: \\

{\it
Supposons que $S$, la partie lin\'eaire de $X$, v\'erifie la condition diophantienne $(\omega)$ de Bruno. Supposons en outre 
que $X$ admette une forme normale de la forme $\hat a. S$, o\`u $a\in \Bbb C[[x_1,\ldots,c_n]]$. Alors, $X$ admet un diff\'eomorphisme 
normalisant holomorphe au voisinage de l'origine dans $\Bbb C^n$.\\
}

D'autre part, J. Vey demontra le r\'esultat suivant \cite{vey-ham}~:\\

{\it
Soient $X_1,\ldots, X_n$, $n$ champs de vecteurs holomorphes au voisinage du point singulier commun $0\in \Bbb C^{2n}$. 
On suppose qu'ils commutent entre eux (i.e. $[X_i,X_j]=0$), qu'ils sont hamiltoniens et que leur partie lin\'eaire 
\`a l'origine sont lin\'eairement ind\'ependantes. Alors, ces champs de vecteurs sont simultan\'ement et holomorphiquement 
normalisables.\\
}

Dans un pr\'ec\'edent travail, nous avons montr\'e que ces r\'esultats ne sont que des faces diff\'erentes d'un m\^eme 
r\'esultat dont nous allons rappeler un des aspects. La premi\`ere id\'ee consiste \`a travailler non pas avec un seul champ 
de vecteurs mais avec une famille de champs de vecteurs commutants entre eux. 
Nous proc\'edons comme suit et nous en profitons pour faire quelques rappels et fixer quelques notations.

Soit $n\geq 2$ un entier et $\lie g$ une alg\`ebre de Lie commutative de dimension finie $l$ sur $\Bbb C$. 
Soient $\lambda_1,\ldots,\lambda_n$ des formes lin\'eaires complexes sur $\lie g$ tel que le morphisme de Lie $S$ de $\lie g$ 
dans l'alg\`ebre  de Lie des champs de vecteurs lin\'eaires de $\Bbb C^n$ d\'efinit par 
$S(g)=\sum_{i=1}^n\lambda_i(g)x_i\frac{\partial}{\partial x_i}$ 
soit injectif. 
Pour tout $Q\in \Bbb N^n$ et $1\leq i\leq n$, on d\'efinit le poids $\alpha_{Q,i}(S)$ de $S$ comme \'etant la forme lin\'eaire 
$\sum_{j=1}^n{q_j\lambda_j(g)}-\lambda_i(g)$. Soit $\|.\|$ une norme sur $\lie g^*$, le $\Bbb C$-espace vectoriel des formes lin\'eaires 
sur $\lie g$.
On d\'efinit la suite de r\'eels positifs 
$$
\omega_k(S)=\inf\{\|\alpha_{Q,i}\|\neq 0, 1\leq i\leq n, 2\leq |Q|\leq 2^k\},
$$
et on dit que $S$ est {\bf diophantien} si la condition suivante est satisfaite~: 
$$
(\omega(S))\quad\quad -\sum_{k\geq 0}\frac{\ln \omega_k(S)}{2^k}<+\infty
$$
Soit $\vfo n k$ (resp. $\fvfo n k$) l'alg\`ebre de Lie des germes de champs de vecteurs holomorphes (resp. formels) d'ordre 
$\geq k$ en $0\in \Bbb C^n$. Soit $\left(\fvfo n 1\right)^S$ (resp. $\left(\widehat{\cal O}_n\right)^S$) 
le commutateur formel  $S$ (resp. l'anneau des integrales premi\`eres formelles), c'est-\`a-dire l'ensemble des champs de vecteurs 
formels $X$ (resp. s\'eries formelles $f$) tel que $[S(g),X]=0$ (resp. ${\cal L}_{S(g)}(f)=0$) pour tout $g\in \lie g$. 

Une d\'eformation nonlin\'eaire $S+\epsilon$ de $S$ est un morphisme de Lie de $\lie g$ dans $\vfo n 1$ tel que 
$\epsilon\in \hom C(\lie g, \vfo n 2)$. 
Soit $\hat\Phi$ un diffeomorphisme formel de $(\Bbb C^n,0)$ que l'on suppose \^etre tangent \`a l'identit\'e en $0$. 
On d\'efinit $\hat \Phi^*(S+\epsilon)(g)\=def\hat\Phi^*(S(g)+\epsilon(g))$ le conjugu\'e de $S+\epsilon$ par $\hat \Phi$.
Dans notre pr\'ec\'edent travail nous avons d\'efini la notion de forme normale formelle de $S+\epsilon$ relativement \`a $S$; 
elle est unique modulo un groupe de transformations formelles. 
Un des r\'esultats principaux que nous avons obtenu est le suivant~:
\begin{theo}\cite{Stolo-sintg2, Stolo-intg-cras}\label{stolo-sintg2}
Soit $S$ un morphisme de Lie diagonal et injectif tel que la condition $(\omega(S))$ soit satisfaite. 
Soit $S+\epsilon$ une d\'eformation non-lin\'eaire holomorphe de $S$. On suppose qu'elle admet  
un \'el\'ement de $\hom C\left(\lie g, \left(\widehat{\cal O}_n\right)^S\otimes_{\Bbb C}S(\lie g)\right)$ comme forme normale formelle. 
Alors elle admet diff\'eomorphisme normalisant holomorphe au voisinage de $0\in\Bbb C^n$.
\end{theo}

Ce r\'esultat admet le r\'esultat de J. Vey et celui de A.D. Bruno comme corollaire. On montre, en particulier, 
que la condition diophantion $\omega(S)$, qui est pr\'ecis\'ement la condtion $(\omega)$ de Bruno dans le cas d'un seul champ de vecteur, 
est automatiquement satisfaite dans le cadre du th\'eor\`eme de Vey. D'ailleurs son \'enonc\'e ne mentionne aucune condition de petits 
diviseurs. La condition diopantienne $(\omega(S))$ est en g\'en\'eral bien plus faible que celle de Bruno.

Dans cet article, nous montrerons qu'il n'est pas n\'ec\'essaire que les champs de vecteurs en question aient tous une 
partie lin\'eaire non-nulle.

\subsection{R\'esultat principal}

Soit $X_1$ un champ de vecteurs holomorphe au voisinage de $0\in \Bbb C^n$, de partie lin\'eaire $s$. Nous dirons 
que $X_1$ est r\'egulier relativement \`a $S$ si~: 
\begin{itemize}
\item $s\in S(\lie g)$ et $\left(\fvfo n 1\right)^S=\left(\fvfo n 1\right)^{s}$;
\item la forme normale formelle de $X_1$ appartient \`a $\left(\widehat{\cal O}_n\right)^S\otimes_{\Bbb C}S(\lie g)$.
\end{itemize}
Comme nous le verrons plus loin, si $X_2,\ldots, X_p$ sont des champs de vecteurs commutants avec $X_1$, alors 
il existe un syst\`eme de coordonn\'ees formelles (unique modulo un groupe de transformations formelles) dans lequel 
les $X_i$ commutent avec $s$. Nous dirons alors que la famille $\{X_1,\ldots,X_l\}$ est {\bf normalis\'ee 
relativement} \`a $s$. La partie junior d'un champ en $0$ est sa partie homog\`ene de plus bas degr\'e.
Nous nous proposons de d\'emontrer le r\'esultat suivant~:
\begin{theo}\label{th-princ}
Soit $S$ un morphisme de Lie diagonal et injectif tel que la condition $(\omega(S))$ soit satisfaite. 
Soit $X_1$ un champs de vecteurs holomorphe au voisinage de son point singulier $0\in \Bbb C^n$. 
On le suppose r\'egulier par rapport \`a $S$. Soient $X_2,\ldots X_l$ des champs de vecteurs 
holomorphes au voisinage de $0$ et commutants avec $X_1$.
On suppose que la famille $\{X_1,\ldots, X_l\}$ admet une forme normale relativement \`a $s$ 
constitu\'ee d'\'el\'ements du $\left(\widehat{\cal O}_n\right)^S$-module 
engendr\'e par $S(\lie g)$ et que leur partie junior forment une famille libre sur $\left(\widehat{\cal O}_n\right)^S$. 
Alors la famille admet un diff\'eomorphisme normalisant holomorphe au voisinage de $0\in\Bbb C^n$.
\end{theo}
Nous appelerons {\bf alg\`ebre de type Cartan} une telle collection de champs de vecteurs, d\'enomination que nous expliquerons 
dans la suite. Nous montrerons qu'un des travaux de H. Ito, concernant les champs hamiltoniens, 
est un cas particulier de notre r\'esultat. 
\begin{rem}
Contrairement \`a la situation du th\'eor\`eme \ref{stolo-sintg2}, 
dans la famille de champs de vecteurs commutants 
$\{X_1,\ldots,X_l\}$, seul $X_1$ doit avoir une partie lin\'eaire non-nulle.
\end{rem}
\begin{rem}
La partie lin\'eaire du champ r\'egulier peut \^etre liouvillienne bien que l'alg\`ebre $S(\lie g)$ 
soit diophantienne. Par exemple, pour $l=2$, $n=4$, $S(g_i)=x_i\frac{\partial }{\partial x_i} - y_i\frac{\partial }{\partial y_i}$, 
$i=1,2$. Le morphisme $S$ est diophantien. Soient $\zeta<0$ un nombre irrationnel liouvillien et 
$s=S(g_1)+\zeta S(g_2)$; ce champ lin\'eaire est r\'egulier par rapport \`a $S$ et 
admet $1,-1,\zeta, -\zeta$ comme valeurs propres. 
On peut donc choisir $\zeta$ de sorte que $s$ ne v\'erifie pas la condition diophantienne de Bruno. 
Ainsi, le caract\`ere diophantien du probl\`eme ne se lit pas directement 
sur la partie lin\'eaire de l'\'el\'ement r\'egulier.
\end{rem}
\begin{rem}
Nous avons appel\'e ces alg\`ebres de champs de vecteurs des "alg\`ebres de type Cartan" par analogie 
avec les sous-alg\`ebres de Cartan d'alg\`ebres de Lie de dimension finie. Une alg\`ebre de type Cartan 
n'est autre qu'un commutateur d'un \'el\'ement r\'egulier, les \'el\'ements de cette alg\`ebre \'etant semi-simples 
sur l'anneau $\widehat{\cal O}^S_n$.
D'autre part, un th\'eor\`eme fondamental stipule que toutes les sous-alg\`ebres de Cartan d'une m\^eme 
alg\`ebre de Lie sont conjugu\'ees entre elles; dans notre situation, les alg\`ebres de type Cartan 
qui sont formellement conjugu\'ees le sont holomorphiquement.
\end{rem}
Ce r\'esultat a \'et\'e anonc\'e dans \cite{Stolo-cartan-cras}.

\section{Notations}

Soit $R=(r_1,\ldots,r_n)\in \left(\Bbb R^*_+\right)^n$; le polydisque ouvert centr\'e en $0\in \Bbb C^n$
 et de polyrayon $R$ sera not\'e $D_R=\{z\in \Bbb C^n\;|\;|z_i|<r_i\}$. Si $r>0$ alors $D_r$ d\'esignera 
le polydisque $D_{(r,\ldots,r)}$. La fronti\`ere distingu\'ee de $D_R$ sera not\'ee ${\cal C}_R$; c'est 
l'ensemble ${\cal C}_R=\{z\in \Bbb C^n\;|\;\forall\, 1\leq i\leq n,\quad |z_i|=R_i\}$. 
Soient $R=(r_1,\ldots,r_n)\in \left(\Bbb R^*_+\right)^n$, $R'=(r'_1,\ldots,r'_n)\in \left(\Bbb R^*_+\right)^n$; 
on \'ecrira $R\leq R'$ lorsque $\forall 1\leq i\leq n$, $r_i\leq r'_i$. Si $Q=(q_1,\ldots,q_n)\in \Bbb Z^n$, 
on notera $|Q|=q_1+\cdots+q_n$; cette application, restreinte \`a $\Bbb N^n$, sera appel\'ee norme.

Si $f\in {\cal O}(\overline{D_R})$ est une fonction holomorphe au voisinage du 
polydisque ferm\'e $\overline{D_R}$, nous poserons $\|f\|_R=\sup_{x\in \overline{D_R}}|f(x)|$.

\subsection{Normes}

Soit $f$ un \'el\'ement de $\Bbb C[[x_1,\ldots,x_n,\frac{1}{x_1},\ldots,\frac{1}{x_1}]]$ que l'on \'ecrira 
$f=\sum_{Q\in \Bbb Z^n}{f_Qx^Q}$; on pose alors $\bar f\=def \sum_{Q\in \Bbb Z^n}{|f_Q|x^q}$. 
L'ordre d'un tel \'el\'ement est le plus petit entier relatif $k\in \Bbb Z$ tel que 
$f_Q\neq 0$ pour $Q\in \Bbb Z^n$ de norme \'egal \`a $k$.
On dira qu'un tel \'el\'ement $g$ domine $f$, si $\forall Q\in \Bbb Z^n,\;|f_Q|\leq |g_Q|$;
dans ce cas, on \'ecrira $f\prec g$. 
Soit $R=(r_1,\ldots,r_n)\in \left(\Bbb R^*_+\right)^n$; on pose 
$$
|f|_R\=def\sum_{Q\in \Bbb Z^n}{|f_Q|R^{Q}}=\bar f(r_1,\ldots,r_n).
$$
On alors les propri\'et\'es suivantes~:
\begin{eqnarray*}
\overline{fg} & \prec & \bar f\bar g\\
\text{if}\;f\prec g & \text{alors} & |f|_{R}\leq |g|_R\\
\overline{\frac{\partial f}{\partial x_k}} & = & \frac{\partial \bar f}{\partial x_k}
\quad \text{lorsque} \quad f\in \Bbb C[[x_1,\ldots,x_n]] \\
\end{eqnarray*}

Soient $v=(v_1,\ldots, v_n)\in \left(\Bbb R^*_+\right)^n$, $0<r<r'$, $R=r.v$, $R'=r'.v$ et $f=\sum_{Q\in \Bbb N^n}{f_Qx^Q}$ une fonction holomorphe au voisinage de 
$\overline{D_{R'}}$, 
on a alors~:
\begin{eqnarray}
\|f\|_{R'} & \leq & |f|_{R'}\\
|f|_R & \leq & \left(\frac{R}{r}\right)^m |f|_{R'}\quad\text{si}\quad \text{ord}(f)\geq m \label{norm-rayon}\\
\left|\frac{\partial f}{\partial z_i}\right|_R & \leq & \frac{d}{rv_i}|f|_R\quad\text{si $f$ est un polyn\^ome de degr\'e}\leq d 
\label{norm-diff}
\end{eqnarray}
\begin{lemm}\cite{Stolo-sintg2}\label{lemm-inverse}
Soient $f,g$ des fonctions holomorphes au voisinage ${\cal U}$ de la fonti\`ere distingu\'ee de 
$D_R$. On suppose que $f$ ne s'annule pas sur ${\cal U}$ et que $\left|\frac{1}{f}\right|_R|g|_R<1$.
Alors,
$$
\left|\frac{1}{f+g}\right|_R\leq 
\left|\frac{1}{f}\right|_R\frac{1}{1-\left|\frac{1}{f}\right|_R|g|_R}
$$ 
\end{lemm}


On posera 
$$
{\cal H}_n(R)=\{f\in \Bbb C[[x_1,\ldots, x_n]]\;|\;|f|_R<+\infty\}
$$
\subsection{Forme normale}

On rappelle le th\'eor\`eme de Poincar\'e-Dulac \cite{Arn2}
\begin{theo}
Soient $X=S+R$ un champs de vecteurs s'annulant \`a l'origine; $S$ sa partie lin\'eaire suppos\'ee semi-simple et $R$ 
un champ de vecteurs non-lin\'eaire. Il existe un diff\'eomorphisme formel $\hat\Phi$ de $(\Bbb C^n,0)$, tangent 
\`a l'identit\'e en $0$ tel que $\hat\Phi^*X=S+N$ o\`u $N$ est un champ de vecteurs formel tel que $[S,N]=0$. 
\end{theo}
On dit alors que $\hat\Phi^*X$ est une forme normale formelle de $X$.
On notera~:
\begin{itemize}
\item $\hvf n m$ l'espace des champs de vecteurs de $\Bbb C^n$ homog\`enes de degr\'e $m$;
\item $\pvf n m {m'}$ l'espace des champs de vecteurs de $\Bbb C^n$ polynomiaux d'ordre $\geq m$ et de degr\'e $\leq m'$;
\item $\vfo n k$ (resp. $\fvfo n k$) l'espace des germes en $0\in \Bbb C^n$ de champs de vecteurs holomorphes 
(resp. formels) d'ordre $\geq k$ en $0$;
\item ${\cal O}_n$ (resp. $\widehat{\cal O}_n$) l'anneau des germes en $0\in \Bbb C^n$ de fonctions holomorphes (resp. formelles);
\item si $X\in \fvfo n 1$ et $k\in \Bbb N^*$, $J^k(X)$ d\'esigne le polyn\^ome de Taylor en $0$ degr\'e $\leq k$.
\end{itemize}

\section{Alg\`ebres de type Cartan}
\subsection{D\'efinitions}

Soit $\lie g$ une $\Bbb C$-alg\`ebre de Lie commutative de dimension $l$. Soit 
$S:\lie g\rightarrow \hvf n 1$ un morphisme de Lie de l'alg\`ebre de Lie $\lie g$ dans l'alg\`ebre 
de Lie $\hvf n 1$ des champs de vecteurs lin\'eaires de $\Bbb C^n$. On supposera que $S$ est injective et
semi-simple; on peut donc supposer qu'il existe des formes lin\'eaires $\lambda_1,\ldots, \lambda_n\in \lie g^*$ telles que 
$$
\forall g\in \lie g, \;S(g)=\sum_{i=1}^n\lambda_i(g)x_i\frac{\partial}{\partial x_i}.
$$
Si $\{g_1,\ldots, g_l\}$ d\'esigna une base de $\lie g$,  on notera $S_i\=def S(g_i)$. On prendra alors norme 
sur $\lie g^*$ d\'efinie par $\|\alpha\|=\max_i|\alpha(g_i)|$.
Le morphisme $S$ d\'efinit une repr\'esentation $\rho_k$ de $\lie g$ dans $\hvf n k$ par $\rho_k(g)p=[S(g),p]$ ($k\geq 2$).
Pour tout $Q\in \Bbb N^n$, $|Q|=k$, $1\leq i\leq n$, les formes lin\'eaires $\alpha_{Q,i}\=def (Q,\lambda)-\lambda_i$ sont 
les poids de cette repr\'esentation.
On a une d\'ecomposition de $\hvf n k$ en sommes directes de espaces de poids de cette repr\'esentation.
$$
\hvf n k=\bigoplus_{\alpha}\left(\hvf n k\right)_{\alpha}(S)
$$
o\`u $\left(\hvf n k\right)_{\alpha(S)}=\{p\in \hvf n k\,|\,\forall g\in \lie g,\,[S(g),p]=\alpha(g)p\}\neq \{0\}$. 
On notera 
\begin{eqnarray*}
\left(\fvfo n 2\right)^S & \=def &\large\{Y\in \vfo n 1\;|\;\forall g\in\lie g,\;[S(g),Y]=0\}=\left(\vfo n 1\right)_0(S),\\
\widehat{\cal O}_n^S & \=def & \{f\in \widehat{\cal O}_n\;|\;\forall g\in \lie g,\; {\cal L}_{S(g)}(f)=0\}.
\end{eqnarray*}
Nous renvoyons le lecteur au chapitre de 5, p.22-29 de \cite{Stolo-sintg2}, pour de plus amples d\'etails.
\begin{defi}
Un \'el\'ement $g_0$ de $\lie g$ est dit {\bf r\'egulier} relativement \`a $S$, si $\left(\fvfo n 1\right)^S=\left(\fvfo n 1\right)^{S(g_0)}$.
\end{defi}

Je dois l'\'enonc\'e suivant \`a Marc Chaperon.
\begin{lemm}
$\lie g$ admet un \'el\'ement r\'egulier relativement \`a $S$. 
\end{lemm}
\begin{Preuve}
Condis\'erons l'ensemble ${\cal P}^*$ des poids non-nuls de $S$ dans $\fvfo n 2$. Il est constitu\'e d'une infinit\'e d\'enombrable 
de forme lin\'eaire non nulle de $\lie g$. Chacune de ces formes lin\'eaires d\'efinit un hyperplan lin\'eaire de $\lie g$.
Ce dernier \'etant un espace vectoriel de dimension finie, il ne saurait \^etre \'egal \`a l'union de l'infinit\'e d\'enombrable 
d'hyperplans d\'efinis par ${\cal P}^*$. Il existe donc un \'el\'ement $g_0\in \lie g$, non nul, qui n'appartient pas \`a cette 
union d'hyperplans. On a donc $\left(\fvfo n 1\right)^{S(g_0)} \hookrightarrow \left(\fvfo n 1\right)^S$. D'o\`u le r\'esultat.
\end{Preuve}

\begin{defi}
Un champ de vecteur holomorphe $X\in \vfo n 1$ sera dit {\bf r\'egulier} relativement \`a $S$ si sa patie lin\'eaire $J^1(X)$ en 
$0$ est r\'eguli\`ere relativement \`a $S$ et si une de ces formes normales formelles appartient au $\widehat{\cal O}_n^S$-module 
engendr\'e par $S(\lie g)$.
\end{defi}

\begin{defi}
Soit $X\in \fvfo n 1$ un champ de vecteur formel de $\Bbb C^n$. On dira que $X$ est normalis\'e (resp. \`a l'ordre $k\geq 2$) 
relativement \`a $S$ s'il appartient \`a $\left(\fvfo n 1\right)^S$ (resp. $\left(\fvfo n 1\right)^S$ mod $\fvfo n {k+1}$).
\end{defi}
\begin{lemm}\label{normal}
Soit $X$ un \'el\'ement r\'egulier relativement \`a $S$. Soit $Y$ un champ de vecteur commutant avec $X$ et s'annulant en $0$.
Si $X$ est normalis\'e \`a l'ordre $k$ alors $Y$ est normalis\'e \`a l'ordre $\text{Ord}(Y)+k -1$.
\end{lemm}
\Preuve
Soit $s$ la partie lin\'eaire de $X_1$. Puisque $X_1$ et $Y$ commutent, chaque composante homog\`ene $l_k$ de degr\'e $k$ 
de $[X_1, Y]$ est nulle. Notons $r$ l'ordre de $Y$ en $0$. On a, pout tout $j\in \Bbb N$, 
$$
l_{r+j}=\sum_{p=1}^{j+1} [X_p,Y_{r+j-p+1}]
$$
o\`u $X_p$ d\'esigne la partie homog\`ene de degr\'e $p$ de $X$. En particulier $0=l_r=[S,Y_r]$.
On montre le r\'esultat par r\'ecurrence sur $k-1\geq j\geq 0$. Pour $j=0$, c'est la remarque pr\'ec\'edente.
Supposons que le r\'esultat est vrai \`a l'ordre $j-1$.
Gr\^ace \`a l'identit\'e de Jacobi, on a 
\begin{eqnarray*}
0=[s,l_{r+j}] & = & \sum_{p=1}^{j+1} [s,[X_p,Y_{r+j-p+1}]]\\
&=& \sum_{p=1}^{j+1} -[X_p,[Y_{r+j-p+1},s]]-[Y_{r+j-p+1},[s,X_p]]
\end{eqnarray*}
Par hypoth\`ese de r\'ecurrence, on a alors 
$$
0=[s,l_{r+j}]=[s,[Y_{r+j},s]];
$$
on obtient le r\'esultat car $[s,.]$ est inversible sur son image (on rappel que $s$ est diagonal).

\cqfd
\begin{coro}
Il existe un diff\'eomorphisme formel $\hat\phi$ tel que $[\hat\phi^*X,s]=[\hat\phi^*Y,s]=0$.
\end{coro}
\begin{defi}
Soit $X_1$ un \'el\'ement r\'egulier relativement \`a $S$. Soient $X_2,\ldots, X_l$ une famille de champs commutants avec 
$X_1$. On dira que la famille $\{X_1,\ldots, X_l\}$ est une {\bf alg\`ebre de type Cartan} relativement \`a $S$, 
si leur forme normale formelle appartiennent au $\widehat{\cal O}_n^S$-module engendr\'e par $S(\lie g)$ et 
si leur patie junior forment famille libre sur $\widehat{\cal O}_n^S$.
\end{defi}

On dira que l'alg\`ebre de type Cartan est normalis\'ee \`a l'ordre $k$ si son \'el\'ement r\'egulier l'est.

\begin{lemm}\label{libre}
Soient $X_1$ un \'el\'ement r\'egulier relativement \`a $S$ et $X_2,\ldots, X_l$ des champs holomorphes commutants avec 
$X_1$. Soit $\hat \Phi$ un diff\'eomorphisme normalisant de $X_1$. Supposons que les $\hat \Phi^*X_i$ appartiennent 
au $\widehat{\cal O}_n^S$-module engendr\'e par $S(\lie g)$ et qu'ils sont libre sur $\widehat{\cal O}_n^S$. 
Posons $\hat \Phi^*X_i=\sum_{i=1}^la_{i,j}S(g_j)$ avec $a_{i,j}\in \widehat{\cal O}_n^S$ et 
$A=(a_{i,j})_{1\leq i,j\leq l}$. Alors $\det A\not\equiv 0$.
\end{lemm}
\Preuve
L'\'equation $\sum_{i=1}^lb_i\hat \Phi^*X_i=0$ avec $b_i\in \widehat{\cal O}_n^S$ est \'equivalente \`a l'ensemble 
d'\'equations $\sum_{i=1}^lb_ia_{i,j}=0$, $j=1,\ldots, l$; ce que l'on peut \'ecrire $A^tB=0$ avec 
$A^t=(a_{j,i})_{1\leq i,j\leq l}$ et $B=(b_i)_{1\leq i\leq l}$. L'application $\widehat{\cal O}_n^S$-lin\'eaire $A^t$ 
du $\widehat{\cal O}_n^S$-module libre $\left(\widehat{\cal O}_n^S\right)^l$ dans lui-m\^eme est injectif si et seulement si 
son determinant $\det A^t=\det A$ n'est pas un diviseur de z\'ero dans $\widehat{\cal O}_n^S$. Or ce dernier est 
un sous anneau de l'anneau $\widehat{\cal O}_n$ qui est int\`egre. Par cons\'equent, la famille $\{\hat \Phi^*X_i\}_{i=1,\ldots, l}$ 
est libre sur $\widehat{\cal O}_n$ si et seulement si $\det A\not\equiv O$ dans $\widehat{\cal O}_n$.
\cqfd
\begin{rem}
En particulier, si $\{X_1,\ldots, X_l\}$ est une alg\`ebre de type Cartan, alors $ord(\det A)=\sum_{i=1}^l(ord X_i-1)$.
\end{rem}
\subsection{Remarques sur la d\'enomination}

La d\'enomination que nous avons adopt\'es repose essentiellement sur une analogie.
Soit $\lie g$ une alg\`ebre de Lie complexe de dimension finie. Un \'el\'ement $g_0\in \lie g$ est r\'egulier si 
la dimension du sous-espace caract\'eristique associ\'e \`a la valeur propre $0$ de $ad_{g_0}$ est minimale dans $\lie g$. 
Une sous-alg\`ebre de Cartan de $\lie g$ est l'ensemble $\{g\in \lie g\;|\; \exists n\in \Bbb N^*,\; ad_{g_0}^n(g)=0\}$ 
o\`u $g_0$ est un \'el\'ement r\'egulier de $\lie g$. En particulier, si $ad_{g_0}$ est semi-simple, l'alg\`ebre de Cartan 
associ\'e n'est autre que le commutateur de $g_0$. D'autre part, un th\'eor\`eme important stipule que les 
sous-alg\`ebres de Cartan sont conjugu\'ees entre elles.

Dans notre situation, $g_0$ est dit r\'egulier relativement \`a $S$ si $S(g_0)$ a le m\^eme commutateur formel que $S$; 
{\it a priori} le commuteur d'un \'el\'ement quelconque de $S(\lie g)$ (qui est semi-simple) le contient. 
Soit $X_1$ un \'el\'ement r\'egulier non-lin\'eaire relativement \`a $S$. Sa forme normale $\hat X_1$ appartient 
$\widehat{\cal O}_n^S$-module engendr\'e par $S(\lie g)$. Elle peut-\^etre vu comme un \'el\'ement semi-simple 
sur $\widehat{\cal O}_n^S$. Soient $X_1,\ldots, X_l$ une alg\`ebre de type Cartan relativement \`a $S$. 
Alors, leurs formes normales formelles commutent \`a $\hat X_1$. Comme corrolaire de notre r\'esultat, 
nous pouvons dire que si deux alg\`ebres de type Cartan sont formellement conjugu\'ees alors elles le sont holomorphiquement.

Il est tr\`es peu problable que l'on puisse d\'efinir, pour les champs de vecteurs singuliers, une notion 
d'alg\`ebres de Cartan ne portant que sur des d\'efinition alg\'ebriques aussi simplement que pour celle d'alg\`ebre de Lie standard); 
en effet, il faut tenir compte, entre autres, des probl\`emes de petits diviseurs.

\section{M\'ethode de Newton diff\'erenti\'ee}
Soientt $X_1$ un \'el\'ement r\'egulier relativement \`a $S$ et $s$ sa partie lin\'eaire.
Soient $\{X_1,\ldots,X_l\}$ une alg\`ebre de type Cartan relativement \`a $S$ et 
$\{S_1,\ldots,S_l\}$ une base de $S(\lie g)$; on note $d_i=Ord(X_i)$; on a $d_1=1$. Dans un bon syst\`eme 
de coordonn\'ees formelles $\hat\Phi$, on a $\hat\Phi^*X_i=\sum_{i=1}^l\hat a_{i,j}S_j$ o\`u $\hat a_{i,j}\in \widehat{\cal O}_n^S$. 
Par hypoth\`eses, leurs parties junior sont libres sur $\widehat{\cal O}_n^S$.
On pose $\hat A=(\hat a_{i,j})_{1\leq i,j\leq l}$ et pour tout entier $p>0$, on pose $A_p=\left(J^{p+d_i-2}(\hat a_{i,j})\right)_{1\leq i,j\leq l}$; 
Par le lemme \ref{libre}, on a $\det A_1\not\equiv 0$.

Supposons l'alg\`ebre de type Cartan normalis\'ee 
\`a l'ordre $m$. On peut donc \'ecrire, pour $i=1,\ldots, l$, $X_i=NF_i^{m+d_i-1}+R_i^{m+d_i}$ o\`u 
$NF_i^{m+d_i-1}$ est un champ de vecteur polynomial, d'ordre $d_i$ et de degr\'e inf\'erieur ou \'egal \`a $m+d_i-1$, 
v\'erifiant $[s,NF_i^{m+d_i-1}]=0$; $R_i^{m+d_i}$ est un champ holomorphe d'ordre $m+d_i$. Par hypoth\`ese, 
on a, pour $i=1,\ldots, l$, 
\begin{equation}
NF_i^{m+d_i-1}=\sum_{j=1}^l{a_{i,j}S_j}\label{forme-normale}
\end{equation}
o\`u $a_{i,j}\=def J^{m+d_i-2}(\hat a_{i,j})\in {\cal O}_n^S$ est un polyn\^ome de degr\'e inf\'erieur 
ou \'egal \`a $m+d_i-2$. Pour simplifier, on posera $A\=def (a_{i,j})_{1\leq i,j\leq l}=A_m\in {\cal M}_l({\cal O}^S)$.

On suppose que $m$ est suffisament grand pour que $\det A\not\equiv 0$ (il suffit donc que $m\geq p_0$). 
Soit $U\in \pvf n {m+1} {2m}$ tel que 
$\exp U^*X_1$ soit normalis\'e \`a l'ordre $2m$. D'apr\`es le lemme \ref{normal}, le $2m+d_i-1$-i\`eme-jet de 
$\exp U^*X_i$ commute avec $s$.
Puisque l'on a 
$$
\exp U^*X_i=X_i+[U,X_i]+\frac{1}{2}[U,[U,X_i]]+\cdots,
$$
il vient que le champ polynomial $J^{2m+d_i-1}\left(R_i^{m+d_i}+[U, NF^{m+d_i-1}_i]\right)$ commute avec $s$.
Rappelons que l'on a une decomposition de Fitting de $\pvf n {m+d_i} {2m+d_i-1}$~:
$$
\pvf n {m+d_i} {2m+d_i-1}=\left(\pvf n {m+d_i} {2m+d_i-1}\right)^+(S)\bigoplus \left(\pvf n {m+d_i} {2m+d_i-1}\right)_0(S)
$$
o\`u $\left(\pvf n {m+d_i} {2m+d_i-1}\right)^+(S)$ (resp. $\left(\pvf n {m+d_i} {2m+d_i-1}\right)_0(S)$) d\'esigne 
la somme directe des espaces de poids non nuls (resp. l'espace de poids nul) de $S$ dans $\pvf n {m+d_i} {2m+d_i-1}$.
Par d\'efinition d'un \'el\'ement r\'egulier, on a $\left(\pvf n {m+d_i} {2m+d_i-1}\right)_0(S)= \left(\pvf n {m+d_i} {2m+d_i-1}\right)_0(s)$.
On en d\'eduit que la projection de $J^{2m+d_i-1}\left(R_i^{m+d_i}+[U, NF^{m+d_i-1}_i]\right)$ sur 
$\left(\pvf n {m+d_i} {2m+d_i-1}\right)^+(S)$ est nulle.

D\'ecomposons $U$ selon les sous-espaces de poids de $\pvf n {m+1} {2m}$. On peut supposer que sa composante $U_0$ selon l'espace 
de poids nul est nulle; en effet, on a $[U, NF^{m+d_i-1}_i]=0$. Soit $\alpha$ un poids non nul de $S$ dans $\pvf n {m+1} {2m}$, 
$U_{\alpha}$ la composante de $U$ le long du sous-espace associ\'e. Un calcul simple, par la formale de Jacobi, montre que 
$[U_{\alpha}, NF^{m+d_i-1}_i]$ appartient \`a l'espace poids de $S$ dans $\fvfo n {2m+d_i-1}$ associ\'e \`a $\alpha$. 
On en d\'eduit que, 
pour tout poids non nul $\alpha$ de $S$ dans $\pvf n {m+1} {2m}$, on a
\begin{equation}
J^{2m+d_i-1}\left(R_{i,\alpha}^{m+d_i}+[U_{\alpha}, NF^{m+d_i-1}_i]\right)=0;\label{equ-cohom0}
\end{equation}
et si $\beta$ est un poids de $S$ dans $\pvf n {m+d_i} {2m+d_i-1}$ qui n'est pas un poids de $S$ dans 
$\pvf n {m+1} {2m}$, alors $J^{2m+d_i-1}\left(R_{i,\beta}^{m+d_i}\right)=0$. Les \'equations $(\ref{equ-cohom0})$ seront 
appel\'ees {\bf \'equations cohomologiques}.

Par cons\'equent, on a , pour tout poids $\alpha$ non nul de $S$ dans $\pvf n {m+1} {2m}$, 
\begin{equation}
[NF_i^{m+d_i-1}, U_{\alpha}] = F_{i,\alpha}^{m+d_i,2m+d_i-1}+R_{i,\alpha}^{2m+d_i} \label{equ-cohom1}
\end{equation}
o\`u l'on a pos\'e $F_{i,\alpha}^{m+d_i,2m+d_i-1}= J^{2m+d_i-1}\left(R_{i,\alpha}^{m+d_i}\right)$ et 
$$
R_{i,\alpha}^{2m+d_i} = [ NF^{m+d_i-1}_i,U_{\alpha}]-J^{2m+d_i-1}\left([ NF^{m+d_i-1}_i,U_{\alpha}]\right)
$$
est un champ polynomial d'ordre sup\'erieur ou \'egal \`a $m+d_i$. 

R\'e\'ecrivons les \'equations (\ref{equ-cohom1}) en utilisant les expressions (\ref{forme-normale}) des 
formes normales partielles. Il vient
\begin{eqnarray*}
F_{i,\alpha}^{m+d_i,2m+d_i-1}+R_{i,\alpha}^{2m+d_i} & = & [NF_i^{m+d_i-1}, U_{\alpha}]\\
&=& \sum_{j=1}^l{a_{i,j}[S_j, U_{\alpha}]+ U_{\alpha}(a_{i,j})S_j}.
\end{eqnarray*}
Ici, $U_{\alpha}(a_{i,j})$ d\'esigne la d\'eriv\'ee de Lie de $a_{i,j}$ le long de $U_{\alpha}$. 
Par d\'efinition, on a $[S_j,U_{\alpha}]=\alpha(g_j)U_{\alpha}$; on r\'e\'ecrit donc ces \'equations sous la forme 
matricielle suivante~:
$$
\begin{pmatrix}
F_{1\alpha}^{m+d_1,2m+d_1-1}+R_{1,\alpha}^{2m+d_1}\\
\vdots \\
F_{l,\alpha}^{m+d_l,2m+d_l-1}+R_{l,\alpha}^{2m+d_l} \\
\end{pmatrix}
=A(x)\begin{pmatrix} \,\alpha(g_1)U_{\alpha} \\ \vdots \\ \,\alpha(g_l)U_{\alpha} \end{pmatrix}
+ \begin{pmatrix} D_1(U_{\alpha}) \\ \vdots \\ D_l(U_{\alpha}) \end{pmatrix}
$$
o\`u $A=(a_{i,j})_{1\leq i,j\leq l}$ et $D_i$ est l'application ${\cal O}_n$-lin\'eaire definie par 
$D_i :U\in \vfo n 2\mapsto \sum_{j=1}^l{U(a_{i,j})S_j}\in \vfo n 2$. 

On \'ecrira pour simplifier $F_i$ pour $F_{i,\alpha}^{m+d_i,2m+d_i-1}$ et $R_i$ pour $R_{i,\alpha}^{2m+d_i}$. 
La matrice $A$ est formellement inversible et 
elle l'est holomorphiquement sur un ouvert ${\cal V}$ dense de $\Bbb C^n$. Soit $C=(c_{i,j})_{1\leq i,j\leq l}$ la matrice 
transpos\'ee des cofacteurs de $A$; on a $C\in {\cal M}_l({\cal O}_n^S)$. En multipliant l'\'equation pr\'ec\'edente 
par $C$, on obtient~:
$$
\begin{pmatrix}
\tilde F_1+\tilde R_1\\
\vdots \\
\tilde F_l+\tilde R_l \\
\end{pmatrix}
=\det A(x)\begin{pmatrix} \,\alpha(g_1)U_{\alpha} \\ \vdots \\ \,\alpha(g_l)U_{\alpha} \end{pmatrix}
+ \begin{pmatrix} \tilde D_1(U_{\alpha}) \\ \vdots \\ \tilde D_l(U_{\alpha}) \end{pmatrix}
$$
o\`u l'on a not\'e 
$$
\tilde F_i  =  \sum_{p=1}^lc_{i,p}F_p,\quad 
\tilde R_i  =  \sum_{p=1}^lc_{i,p}R_p\; \text{ et}\quad
\tilde D_i(U)  =  \sum_{p=1}^lc_{i,p}D_p(U)
$$

Les op\'erateurs $\tilde D_i$ sont nilpotents; on a $\tilde D_i\circ \tilde D_i=0$. En effet, pour tout $1\leq j\leq l$, on a
\begin{eqnarray*}
\tilde D_j(\tilde D_j(U)) & = & \sum_{k=1}^l{c_{j,k}D_k(\tilde D_j(U))}
 =  \sum_{k=1}^l{c_{j,k}D_k\left(\sum_{i=1}^l{c_{j,i}D_i(U)}\right)}\\
& = & \sum_{i,k=1}^l{c_{j,k}c_{j,i}D_k(D_i(U))}\quad \text{par ${\cal O}_n$-lin\'earit\'e des op\'erateurrs $D_i$}.
\end{eqnarray*}
Montrons alors que, pour tout entiers $1\leq i,k\leq l$, $D_k\circ D_i=0$. En effet, 
\begin{eqnarray*}
D_k(D_i(U)) & = & \sum_{p=1}^l{D_i(U)(a_{k,p})S_p}\\
& = & \sum_{p=1}^l{\sum_{q=1}^l{U_{\alpha}(a_{i,q})S_q(a_{k,p})}S_p}.
\end{eqnarray*}
Or, $a_{k,p}\in {\cal O}^{S}_n$ donc $S_q(a_{k,p})=0$; il s'en suit que $D_k\circ D_i=0$ et donc que 
$\tilde D_j\circ \tilde D_j=0$ comme annonc\'e.

En se pla\c{c}ant dans le corps des fractions de $\widehat{\cal O}_n$, on obtient, pour tout $1\leq i\leq l$, 
\begin{equation}
U_{\alpha}= \left(Id-\frac{1}{\alpha(g_i)\det A(x)}\tilde D_i\right)\frac{\tilde F_i+\tilde R_i}{\alpha(g_i)\det A(x)};
\end{equation}
ce que l'on peut \'ecrire 
\begin{equation}
U_{\alpha}(x) =  \frac{1}{\alpha(g_i)\det\;(A(x))}\left(G_i(x)-\frac{1}{\alpha(g_i)\det\;(A(x))}H_i(x)\right)\label{sol-eq-hom}
\end{equation}
avec 
\begin{eqnarray}
G_i & = & \tilde F_i+\tilde R_i = \sum_{p=1}^lc_{i,p}(F_p+R_p)\label{tilde-Fi}\\
H_i & = & \tilde D_i(\tilde F_i+\tilde R_i) = \sum_{q=1}^l c_{i,q}D_q(\tilde F_i+\tilde R_i)= \sum_{p,q=1}^l c_{i,p}c_{i,q}D_q(F_p+R_p)\nonumber\\
& = & \sum_{p,q,r=1}^l c_{i,p}c_{i,q}(F_p+R_p)(a_{q,r})S_r\label{tilde-Di}
\end{eqnarray}

Le second membre de l'\'equation (\ref{sol-eq-hom}) d\'efinit une application holomorphe sur l'ouvert ${\cal V}$
, qui est dense dans $\Bbb C^n$ et o\`u la matrice $A$ est holomorphiquement inversible. L'application $U_{\alpha}$ 
est holomorphe dans $\Bbb C^n$. L'\'equation (\ref{sol-eq-hom}) a donc un sens dans $\Bbb C^n$ et le second membre 
est un polyn\^ome de degr\'e $\leq 2m$ et d'ordre $\geq m+1$. 

\subsection{Domination de la solution des l'\'equations cohomologiques}

Soit $t$ une variable complexe et $x\in \Bbb C^n$; on pose $t.x=(tx_1,\ldots,tx_n)$. Soit ${\cal M}_n$ le corps de fractions de 
$\Bbb C[x_1,\ldots,x_n]$. Les \'equations (\ref{sol-eq-hom}) induisent, dans $\left({\cal M}_n\right)^n[[t]]$, les \'equations
\begin{equation}
U_{\alpha}(t.x) =  \frac{1}{\alpha(g_i)\det\;(A(t.x))}\left(G_i(t.x)-\frac{1}{\alpha(g_i)\det\;(A(t.x))}H_i(t.x)\right)\label{sol-frac}
\end{equation}
L'ordre d'un \'el\'ement $f=\sum_{p\geq 0}a_pt^p\in \left({\cal M}_n\right)^n[[t]]$ est le plus petit entier $k_0$ tel que $a_{k_0}\neq 0$; 
on le note $\text{ord}_t(f)$ et si $F=(f_1,\ldots,f_l)\in \left({\cal M}_n[[t]]\right)^l$.

\begin{lemm}\label{rest-t}
Avec les notations pr\'ec\'edentes, on a 
\begin{itemize}
\item 
$\text{ord}_t\left(\frac{1}{\det A}\sum_{p=1}^lc_{i,p}R_p(t.x)\right)\geq 2m+1$, 
\item 
$\text{ord}_t\left(\frac{1}{\det A^2}\sum_{p,q,r=1}^l c_{i,p}c_{i,q}R_p(a_{q,r})S_r(t.x)\right)\geq 2m+1$. 
\end{itemize}
\end{lemm}
\Preuve
Soit $D=(d_{i,j})_{1\leq i,j\leq l}$ la matrice de cofacteurs de $A$. Par d\'efinition, $d_{i,j}$ est le d\'eterminant de la matrice obtenu 
\`a partir de $A$ en enlevant la $i$-i\`eme ligne et la $j$-i\`eme colonne. Par hypoth\`ese, le coefficient $a_{i,j}$ de $A$ 
est d'ordre sup\'erieur ou \'egal \`a $d_i-1$. La formule du determinant montre alors clairement que l'ordre de 
$d_{i,j}$ est sup\'erieur ou \'egal \`a $\sum_{p\neq i}(d_p-1)$. On en d\'eduit que l'ordre de $c_{i,j}$, o\`u 
$C=(c_{i,j})_{1\leq i,j\leq l}$ est la transpos\'ee de $D$, est sup\'erieur ou \'egal \`a $\sum_{p\neq j}(d_p-1)$.
On rappelle que $R_j$ est d'ordre sup\'erieur ou \'egal \`a $2m+d_j$; donc, $c_{i,j}R_j$ est d'ordre sup\'erieur ou \'egal 
\`a $2m+1+\sum_{p= 1}^l(d_p-1)$. D'autre part, cette m\^eme formule du determinant montre que 
l'ordre de $\det A$ est \'egal \`a $\sum_{p= 1}^l(d_p-1)$ car les termes juniors sont libres sur ${\cal O}_n^S$. On en d\'eduit que 
$$
\text{ord}_t\left(\frac{1}{\det A}\sum_{p=1}^lc_{i,p}(R_p)(t.x)\right)\geq 2m+1+\sum_{p= 1}^l(d_p-1)-\sum_{p= 1}^l(d_p-1)\geq 2m+1
$$

De m\^eme, l'ordre de $\frac{1}{\det A^2}\sum_{p,q,r=1}^l c_{i,p}c_{i,q}(R_p)(a_{q,r})S_r(t.x)$ est sup\'erieur ou \'egal 
\`a 
$$
\sum_{j\neq p}(d_j-1)+\sum_{k\neq q}(d_k-1)+2m+d_p+d_q-1-2\sum_{p= 1}^l(d_p-1)\geq 2m+1
$$
\cqfd

Puisque $U_{\alpha}(t.x)$ est un polyn\^ome de degr\'e $\leq 2m$, ce lemme montre que $U_{\alpha}$ se calcul uniquement \`a partir des $F_i$; et nous allons donner une estimation 
de $U_{\alpha}$ en fonction de celles des $F_i$.
%
%
Par cons\'equent, 
$$
U_{\alpha}=J^{2m}\left({\frac{1}{\alpha(g_i)\det\;(A(x))}}\left({\tilde F_i(x)}-{\frac{1}{\alpha(g_i)\det\;(A(x))}}{\tilde D_i(\tilde F_i)(x)}\right)\right).
$$
On obtient aors la relation de domination suivante~:
\begin{equation}
U_{\alpha}(x)  \prec   \overline{\frac{1}{\alpha(g_i)\det\;(A(x))}\left(\tilde F_i(x)-\frac{1}{\alpha(g_i)\det\;(A(x))}\tilde D_i(\tilde F_i)(x)\right)}\label{sol-frac3};
\end{equation}

Soit $t_0\in ]0,1]$ et $\kappa=(\kappa_1,\ldots ,\kappa_n)\in \left(\Bbb R^*_+\right)^n$. On pose 
$t_0^{\kappa}.x\=def (t_0^{\kappa_1}x_1,\ldots ,t_0^{\kappa_n}x_n)$. Soient $f\in \Bbb C[[x_1,\ldots,x_n]]$ une s\'erie formelle 
et $g\in \Bbb C[[x_1,\ldots, x_n,\frac{1}{x_1},\ldots,\frac{1}{x_n}]]$ telles que $f\prec g$; alors on a $f(t_0^{\kappa}.x)\prec g(t_0^{\kappa}.x)$.

On peut donc \'ecrire, d'apr\`es la relation de domination (\ref{sol-frac3}), la relation de domination suivante~:
\begin{equation}
U_{\alpha}(t_0^{\kappa}.x) \prec  \overline{\frac{1}{\alpha(g_i)\det\;(A(t_0^{\kappa}.x))}}\left(\overline{\tilde F_i(t_0^{\kappa}.x)}+\overline{\frac{1}{\alpha(g_i)\det\;(A(t_0^{\kappa}.x))}}\overline{\tilde D_i(\tilde F_i)(t_0^{\kappa}.x)}\right)\label{sol-frac4}.
\end{equation}

\subsection{Estimation de la solution des \'equations cohomologiques}

On \'ecrit $A=A_{p_0}+R$ o\`u comme pr\'ec\'edement $A_{p_0}=\left(J^{p_0+d_i-2}(\hat a_{i,j})\right)_{1\leq i,j\leq l}$ 
et $R\=def(r_{i,j})_{1\leq i,j\leq l}$ ($p=0=1$ !). 
On \'ecrira 
$$
\begin{pmatrix}
NF^{p_0}_{1}\\
\vdots \\
NF^{p_0+d_l-1}_{l} \\
\end{pmatrix}
\=def A_{p_0}(x)\begin{pmatrix} \,S_1 \\ \vdots \\ \,S_l \end{pmatrix},
$$
de sorte que pour tout $i=1,\ldots,l$, $NF^{m+d_i-1}_i-NF^{p_0+d_i-1}_{i}$ est d'ordre $\geq p_0+d_i$.

\begin{lemm}\label{estim-det}
Soient $\kappa=(\kappa_1,\ldots,\kappa_n)\in \left(\Bbb R^+_*\right)^l$ un vecteur incommensurable,
c'est-\`a-dire que les $\kappa_i$ sont lin\'eairement ind\'ependants sur $\Bbb Q$ et $r>0$. Il existe alors  
$0<t_0\leq 1$ et $\eta>0$ tels que, si $|R|_{t_0^{\kappa}.r}<\eta$ alors on a
$$
\left|\frac{1}{\det\;(A(x))}\right|_{t_0^{\kappa}.r}\leq \frac{c}{t_0^dr^s}
$$
( on rappelle que $t_0^{\kappa}.r=(t_0^{\kappa_1}r,\ldots,t_0^{\kappa_l}r)$). 
En outre, on a $\left|\det\;(A(x))\right|_{t_0^{\kappa}.r}\leq M_1|A_{p_0}|_{t_0^{\kappa}.r}^l+\frac{1}{2c}r^{s}t_0^d$,
o\`u $M_1$ est une constante universelle.
Les constantes positives $c,d, t_0$ et $\eta$ ainsi que l'entier $s$ 
ne d\'ependent que de $\kappa$, de la matrice $A_{p_0}$ et du morphisme $S$.
\end{lemm}
\Preuve

Soit ${\cal S}_l$ le groupe des permutations de $\{1,\ldots,l\}$. Si $\sigma\in {\cal S}_l$, 
$\epsilon(\sigma)$ d\'esigne la signature de $\sigma$.
On rappelle l'expression du deteminant de $A$~: 
\begin{eqnarray*}
\det (A(x)) & = & \sum_{\sigma\in {\cal S}_l}{\epsilon(\sigma)\left(\prod_{i=1}^l{a_{i,\sigma(i)}(x)}\right)}\\
& = & \det A_{p_0} + P(R)\\
\end{eqnarray*}
o\`u $P(Z)\in {\cal O}_n[Z_1,\ldots,Z_{l^2}]$ est un polyn\^ome de $l^2$ variables sans terme constant 
et de degr\'e $l$. Ces coefficients $d_Q(A_{p_0})$ sont des polyn\^omes de degr\'e inf\'erieur ou \'egal \`a $l-1$ 
en les coefficients de $A_{p_0}$. On peut donc \'ecrire  
$P(Z)=\sum_{\substack{ Q\in \Bbb N^{l^2}\\ 0<|Q|\leq l}}{d_Q(A_{p_0})Z^Q}$ ainsi que 
$\det A_{p_0}=\sum_{\substack{T\in \Bbb N^n\\|S|= N}}{p_Tx^T}$ (on rapelle que $\det A_{1}$ est 
un polyn\^ome homog\`ene).

Soit $\kappa=(\kappa_1,\ldots,\kappa_n)\in \left(\Bbb R^+_*\right)^n$ un nuplet incommensurable. On a alors
$$
\det A_{p_0}(t^{\kappa}.x)=\sum_{\substack{T\in \Bbb N^n\\|T|= N}}{p_Tt^{(\kappa,T)}x^T}
\quad\text{o\`u}\quad(\kappa,T)=\sum_{i=1}^n{\kappa_it_i}\text{ et }T=(t_1,\ldots,t_n).
$$
Posons
$$
d=\min_{\substack{T\in \Bbb N^n\\p_T\neq 0,|T|= N}}{(\kappa,T)}.
$$
Gr\^ace \`a l'incommensurabilit\'e de $\kappa$, ce minimum ne peut \^etre atteint qu'une seul fois sur l'ensemble 
$\{T\in \Bbb N^n,p_T\neq 0\;|\;|T|= N\}$. Soit $T_0$ l'\'el\'ement qui le r\'ealise. 
On en d\'eduit que $\det A_{p_0}(t^{\kappa}.x)=p_{T_0}x^{T_0}t^d+t^dq(x,t)$ o\`u
$q(x,t)$ est un polyn\^ome en $x$ de degr\'e $N$ et dont les coefficients sont des fonctions de la variable $t$ et s'annulant
en $t=0$. Par cons\'equent, il existe $0< t_0\leq 1$ (on peut supposer $t_0$ plus petit que $1$) tel que, 
$|t_0^dq(x,t_0)|_{r}\leq \frac{|p_{T_0}|}{4}r^{|T_0|}t_0^d$.

Soit $\tilde P(Z)\in \Bbb C[Z_1,\ldots,Z_{l^2}]$ le polyn\^ome d\'efinit par 
$\tilde P(z)=\sum_{\substack{ Q\in \Bbb N^{l^2}\\ 0<|Q|\leq l}}{\|d_Q(A_{p_0}(x))\|_1 Z^Q}$.
En tant que fonction de $Z$, on a $P\prec \tilde P$. Puisque $\tilde P(0)=0$, 
il existe $\eta>0$ tel que $|\tilde P(Z)|_{\eta}\leq \frac{|p_{T_0}|}{4}r^{|T_0|}t_0^d$; ici, nous avons \'ecrit 
$|\tilde P(Z)|_{\eta}=\sum_{\substack{ Q\in \Bbb N^{l^2}\\ |Q|\leq l}}{\|d_Q(A_{p_0}(x))\|_1\eta^{|Q|}}$.
Par cons\'equent, si $|R(x)|_{t_0^{\kappa}.r}<\eta$ 
alors, puisque $P(R(t_0^{\kappa}.x))\prec \tilde P(\bar R(t_0^{\kappa}.x))$,
on a $|P(R(x))|_{t_0^{\kappa}.r}\leq \frac{|p_{S_0}|}{4}r^{|S_0|}t_0^d$.

Il s'en suit, d'apr\`es le lemme $\ref{lemm-inverse}$, 
\begin{eqnarray*}
\left|\frac{1}{\det\;(A(t_0^{\kappa}.x))}\right|_{r}& \leq & 
\left|\frac{1}{p_{S_0}x^{S_0}t_0^d}\right|_{r}
\left(\frac{1}{1-\left|\frac{1}{p_{S_0}x^{S_0}t_0^d}\right|_{r}|t_0^dq(x,t_0)+P(R(t_0^{\kappa}.x)|_{r}}\right)\\
&\leq & \frac{1}{|p_{S_0}|r^{|S_0|}t_0^d}
\left(\frac{1}{1-\frac{1}{|p_{S_0}|r^{|S_0|}t_0^d}(|t_0^dq(x,t_0)|_{r}+|P(R(t_0^{\kappa}.x)|_{r})}\right)\\
&\leq & \frac{2}{|p_{S_0}|r^{|S_0|}t_0^d}.
\end{eqnarray*}
De plus, on a  
\begin{eqnarray*}
\left|\det\;(A(t_0^{\kappa}.x))\right|_{r} & \leq & |\det A_{P_0}(t_0^{\kappa}.x)|_{r}+ |P(R(t_0^{\kappa}.x))|_{r}\\
& \leq & M_1|A_{P_0}|_{t_0^{\kappa}.r}^l+\frac{|p_{S_0}|}{4}r^{|S_0|}t_0^d
\end{eqnarray*}
o\`u $M_1$ est une constante universelle (on peut prendre $\sharp {\cal S}_l$).
\cqfd\\
L'id\'ee de consid\'erer des polydisques asym\'etriques est due \`a H. Ito. Le r\'esultat suicant est fondamental pour la suite.
\begin{theo}\label{estim-cohom}
Soit $\kappa=(\kappa_1,\ldots,\kappa_n)\in \left(\Bbb R_+^*\right)^n$ un vecteur incommensurable et $t_0$ comme dans 
le lemme \ref{estim-det} et on pose $m=2^k$.
Il existe des constantes $\eta_1>0$ et $c_1>0$ telles que, 
si $1/2<r\leq 1$, $\max_i |NF^{m+d_i-1}_i-NF^{p_0+d_i-1}_i|_{t_0^{\kappa}.r}<\eta_1$ et 
$\max_i \left|D\left(NF^{m+d_i-1}_i-J^1(NF^{m+d_i-1}_i)\right)\right|_{t_0^{\kappa}.r}< \eta_1$, alors, pour tout poids 
non nul de $S$ dans $\pvf n {m+1} {2m}$, on a
$$
|U_{\alpha}|_{t_0^{\kappa}.r}< \frac{c_1}{\omega_{k+1}(S)^2}\max_i|R_{i,\alpha}^{m+d_i}|_{t_0^{\kappa}.r}.
$$
\end{theo}
\Preuve
En utilisant la relation de domination (\ref{sol-frac4}) ainsi que les propri\'et\'e de la norme $|.|_{t_0^K, r}$, 
sous l'hypoth\`se $|R|_{t_0^{\kappa}.r}<\eta$, 
on obtient gr\^ace au lemme \ref{estim-det}
\begin{eqnarray}
|U_{\alpha}|_{t_0^{\kappa}.r}
& \leq & 
\left|\frac{1}{|\alpha(g_i)|\det\;(A(t_0^{\kappa}.x))}\right|_{r}^2
\left( |\alpha(g_i)|\left|\det\;A\right|_{t_0^{\kappa}.r}\left|\tilde F_i\right|_{t_0^{\kappa}.r}
+\left|\tilde D_i(\tilde F_i)\right|_{t_0^{\kappa}.r}\right)\nonumber\\
& \leq & \frac{c^2}{|\alpha(g_i)|^2t_0^{2d}r^{2s}}\nonumber\\
& & \times
\left( 
|\alpha(g_i)| \left( M_1|A_{p_0}|_{t_0^{\kappa}.r}^l+\frac{|p_{S_0}|}{4}r^{|S_0|}t_0^d\right)\left|\tilde F_i\right|_{t_0^{\kappa}.r}
+\left|\tilde D_i(\tilde F_i)\right|_{t_0^{\kappa}.r}
\right)\label{1st-etim-hom-eq}
\end{eqnarray}

Il nous reste \`a estimer $\left|\tilde F_i\right|_{t_0^{\kappa}.r}$ et 
$\left|\tilde D_i(\tilde F_i)\right|_{t_0^{\kappa}, r}$ en fonction des donnn\'ees; les fonctions $\tilde F_i$ et $\tilde D_i(\tilde F_i)$ \'etant 
d\'efinies respectivement par $(\ref{tilde-Fi})$ et $(\ref{tilde-Di})$. 

On rappelle que chaque fonction $c_{i,j}$ (\'el\'ement de la matrice transpos\'ee des cofacteurs) s'exprime comme polyn\^ome 
universel de degr\'e $l-1$ en les $a_{i,j}$. 
Il y a donc une constante universelle $M>0$ telle que $|c_{i,j}|_{t_0^{\kappa}.r}\leq M |A|_{t_0^{\kappa}.r}^{l-1}$. 
On obtient alors les in\'egalit\'es suivantes~:
\begin{eqnarray*}
|\tilde F_i|_{t_0^{\kappa}.r} & \leq & lM|A|_{t_0^{\kappa}.r}^{l-1}|F|_{t_0^{\kappa}.r}\\
|\tilde D_i(\tilde F_i)|_{t_0^{\kappa}.r} & \leq & l^3M^2|A|_{t_0^{\kappa}.r}^{2(l-1)}\max_{1\leq j\leq l}|S_{j}|_{t_0^{\kappa}.r}\max_{1\leq p,q,r\leq l}|F_{p}(a_{q,r})|_{t_0^{\kappa}.r}\\
\end{eqnarray*}

Or, on a $|F_{p}(a_{q,r})|_{t_0^{\kappa}.r}\leq n|F_{p}|_{t_0^{\kappa}.r}
\max_{1\leq i\leq l}\left|\frac{\partial a_{q,r}}{\partial x_i}\right|_{t_0^{\kappa}.r}$. 
On en d\'eduit la majoration suivante~:
$$
|\tilde D_i(\tilde F_i)|_{t_0^{\kappa}.r} \leq  nl^3M^2|A|_{t_0^{\kappa}.r}^{2(l-1)}\max_{1\leq j\leq l}|S_{j}|_{1}
|D(A)|_{t_0^{\kappa}.r}|F|_{t_0^{\kappa}.r},
$$
o\`u $D(A)$ d\'esigne la d\'eriv\'ee de $A(x)$.

Ainsi, sous l'hypoth\`ese $|R|_{t_0^{\kappa}.r}<\eta$, on a 
\begin{eqnarray*}
|U_{\alpha}|_{t_0^{\kappa}.r} & \leq &
 \frac{c^2lM|A|_{t_0^{\kappa}.r}^{l-1}|F|_{t_0^{\kappa}.r}}{|\alpha(g_i)|^2t_0^{2d}r^{2s}}\\
& & \times \left( 
|\alpha(g_i)| \left( M_1|A_{p_0}|_{t_0^{\kappa}.r}^l+\frac{1}{2c}r^{s}t_0^d\right)
+nl^2M|A|_{t_0^{\kappa}.r}^{l-1}\max_{1\leq j\leq l}|S_{j}|_{1}|D(A)|_{t_0^{\kappa}.r}
\right).
\end{eqnarray*}

Il nous reste \`a majorer les normes de $A$ et $D(A)$ en fonctions de celles des formes normales partielles.
Par d\'efinition,  on a , pour tout entier $1\leq i\leq l$,
\begin{eqnarray*}
S_j & = & \sum_{k=1}^n{\lambda_{j,k}x_k\frac{\partial}{\partial x_k}},\\
NF^{m+d_i-1}_i & = & \sum_{j=1}^l{a_{i,j}S_j}\=def\sum_{k=1}^n{x_k g_{i,k}\frac{\partial}{\partial x_k}}
\quad\text{avec}\quad g_{i,k}= \left(\sum_{j=1}^l{\lambda_{j,k}a_{i,j}}\right)\\
NF^{p_0+d_i-1}_{i} & = & \sum_{j=1}^l{J^{p_0+d_i-2}(a_{i,j})S_j}\=def\sum_{k=1}^n{x_k g_{i,k,p_0}\frac{\partial}{\partial x_k}}
\quad\text{avec}\quad g_{i,k,p_0}= \left(\sum_{j=1}^l{\lambda_{j,k}J^{p_0+d_i-2}(a_{i,j})}\right);\\
\end{eqnarray*}
ce que l'on peut r\'e\'ecrire sous la forme
$$
\begin{pmatrix} g_{i,1}\\ \vdots\\ \vdots\\ g_{i,n} \end{pmatrix}
=
\begin{pmatrix} \lambda_{1,1} & \ldots & \lambda_{l,1}\\ \vdots & & \vdots \\ \vdots & & \vdots \\
\lambda_{1,n} & \ldots & \lambda_{l,n} \end{pmatrix}
\begin{pmatrix} a_{i,1}\\ \vdots\\ a_{i,l} \end{pmatrix}
\quad\text{o\`u}\quad l\leq n, \text{ et} 
$$

$$
\begin{pmatrix} g_{i,1,p_0}\\ \vdots\\ \vdots\\ g_{i,n,p_0} \end{pmatrix}
=
\begin{pmatrix} \lambda_{1,1} & \ldots & \lambda_{l,1}\\ \vdots & & \vdots \\ \vdots & & \vdots \\
\lambda_{1,n} & \ldots & \lambda_{l,n} \end{pmatrix}
\begin{pmatrix} J^{p_0+d_i-1}(a_{i,1})\\ \vdots\\ J^{p_0+d_i-1}(a_{i,l}) \end{pmatrix}
\quad\text{o\`u}\quad l\leq n. 
$$
Puisque les $S_i$ sont lin\'eairement ind\'ependants sur $\Bbb C$, la matrice 
$(\lambda_{j,i})_{\substack{1\leq i\leq n\\ 1\leq j\leq l}}$ est de rang $l$. Sans nuire \`a la g\'en\'eralit\'e, 
on peut supposer que la matrice $L\=def (\lambda_{j,i})_{1\leq i,j\leq l}$ est inversible d'inverse $L^{-1}\=def (\tilde\lambda_{i,j})_{1\leq i,j\leq l}$. 
On peut alors \'ecrire, pour $1\leq i,j\leq l$,
$$
a_{i,j}(x)- J^{p_0+d_i-2}(a_{i,j}) =  \sum_{k=1}^{l}{\tilde \lambda_{j,k}(g_{i,k}(x)-g_{i,k,p_0})}
 \prec  \sum_{k=1}^{l}{|\tilde \lambda_{j,k}| \overline{(g_{i,k}(x)-g_{i,k,p_0})}},
$$
Puisque $1/2<r$, on a
\begin{eqnarray*}
\left|g_{i,k}-g_{i,k,p_0}\right|_{t_0^{\kappa}.r}& \leq & \frac{2}{t_0^{\kappa_k}}(t_0^{\kappa_k}r) \left|g_{i,k}-g_{i,k,p_0}\right|_{t_0^{\kappa}.r}\\
& = & \frac{2}{t_0^{\kappa_k}}\left|x_k(g_{i,k}-g_{i,k,p_0})\right|_{t_0^{\kappa}.r}\\
& \leq & \frac{2}{\min_k t_0^{\kappa_k}}\max_i|NF^{m+d_i-1}_i-NF^{p_0+d_i-1}_i|_{t_0^{\kappa}.r}.
\end{eqnarray*}
Par cons\'equent, on a 
\begin{equation}\label{estim-A}
|A-A_{p_0} |_{t_0^{\kappa}.r}\leq \frac{2l|L^{-1}|}{\min_k t_0^{\kappa_k}}\max_i|NF^{m+d_i-1}_i-NF^{p_0+d_i-1}_i|_{t_0^{\kappa}.r}.
\end{equation}
D'autre part, pour tout $1\leq k\leq n$, on a
$$
\frac{\partial x_k\overline{(g_{i,k}(x)-g_{i,k}(0))}}{\partial x_p}=\delta_{k,p}\overline{(g_{i,k}(x)-g_{i,k}(0))}
+x_k\frac{\partial \overline{(g_{i,k}(x)-g_{i,k}(0))}}{\partial x_p},
$$
c'est-\`a-dire 
$$
x_k\frac{\partial \overline{(g_{i,k}(x)-g_{i,k}(0))}}{\partial x_p}=\frac{\partial x_k\overline{(g_{i,k}(x)-g_{i,k}(0))}}{\partial x_p}
-\delta_{k,p}\overline{(g_{i,k}(x)-g_{i,k}(0))}.
$$
Or $x_k\frac{\partial \overline{(g_{i,k}(x)-g_{i,k}(0))}}{\partial x_p}$, 
$\frac{\partial x_k\overline{(g_{i,k}(x)-g_{i,k}(0))}}{\partial x_p}$ and $\overline{(g_{i,k}(x)-g_{i,k}(0))}$ sont
des series formelles \`a coefficients positifs; par cons\'equent, on a 
$$
x_k\frac{\partial \overline{(g_{i,k}(x)-g_{i,k}(0))}}{\partial x_p}\prec \frac{\partial x_k\overline{(g_{i,k}(x)-g_{i,k}(0))}}{\partial x_p}.
$$

Pusique $1/2<r$, on a 
\begin{eqnarray*}
\left|\frac{\partial \overline{(g_{i,k}(x)-g_{i,k}(0))}}{\partial x_p}\right|_{t_0^{\kappa}.r} & \leq & \frac{2}{t_0^{\kappa_k}}(t_0^{\kappa_k}r) \left|\frac{\partial \overline{(g_{i,k}(x)-g_{i,k}(0))}}{\partial x_p}\right|_{t_0^{\kappa}.r} \\ 
& \leq & \frac{2}{\min_k t_0^{\kappa_k}}\left|x_k\frac{\partial \overline{(g_{i,k}(x)-g_{i,k}(0))}}{\partial x_p}\right|_{t_0^{\kappa}.r} \\
& \leq & \frac{2}{\min_k t_0^{\kappa_k}}\left|\frac{\partial x_k\overline{(g_{i,k}(x)-g_{i,k}(0))}}{\partial x_p}\right|_{t_0^{\kappa}.r}.
\end{eqnarray*}
On en d\'eduit que, pour tout $1\leq i,j\leq l$,
\begin{eqnarray*}
\left|\frac{\partial \overline{(a_{i,j}(x)-a_{i,j}(0))}}{\partial x_p}\right|_{t_0^{\kappa}.r} & \leq & \sum_{k=1}^{l}{|\tilde \lambda_{j,k}| \left|\frac{\partial \overline{(g_{i,k}(x)-g_{i,k}(0))}}{\partial x_p}\right|_{t_0^{\kappa}.r} }\\
& \leq & \frac{2}{\min_k t_0^{\kappa_k}}\sum_{k=1}^{l}{|\tilde \lambda_{j,k}| \left|\frac{\partial x_k\overline{(g_{i,k}(x)-g_{i,k}(0))}}{\partial x_p}\right|_{t_0^{\kappa}.r}}\\
& \leq & \frac{2l|L^{-1}|}{\min_k t_0^{\kappa_k}}\max_i\left|D\left(NF^{m+d_i-1}-J^1(NF^{m+d_i-1})\right)\right|_{t_0^{\kappa}.r};
\end{eqnarray*}
c'est-\`a-dire
\begin{equation}\label{estim-DA}
|D(A)|_{t_0^{\kappa}.r} \leq \frac{2l|L^{-1}|}{\min_k t_0^{\kappa_k}}\max_i\left|D\left(NF^{m+d_i-1}-J^1(NF^{m+d_i-1})\right)\right|_{t_0^{\kappa}.r}.
\end{equation}

Posons $\eta_1=\frac{\min_k t_0^{\kappa_k}}{2l|L^{-1}|}\eta$;
si $\max_i\left|NF^{m+d_i-1}-NF^{p_0+d_i-1}\right|_{t_0^{\kappa}.r}<\eta_1$ et \linebreak 
$\max_i\left|D\left(NF^{m+d_i-1}-J^1(NF^{m+d_i-1})\right)\right|_{t_0^{\kappa}.r}< \eta_1$, 
alors $|A-A_{p_0}|_{t_0^{\kappa}.r}<\eta$ et $|D(A)|_{t_0^{\kappa}.r}\leq \eta$. Posons $C=|A_{p_0}|_{1}$; on a alors 
$|A|_{t_0^{\kappa}.r}\leq C+\eta$. De plus, $\alpha$ \'etant un poids de $S$ d'ordre $\leq 2m=2^{k+1}$, 
soit $1\leq i\leq l$ tel que $|\alpha(g_i)|= \max_j |\alpha(g_j)|=\|\alpha\|$; on a alors 
$\omega_{k+1}(S)\leq |\alpha(g_i)|$. De plus, la pr\'esence de petits diviseurs nous permet de supposer que $\omega_k(S)\leq 1$.

Posons alors 
$$
c_1  =  \frac{2^{2s}c^2lM(C+\eta)^{l-1}}{t_0^{2d}}
\left( 
 \left( M_1C^l+\frac{1}{2c}\right)
+l^2M(C+\eta)^{l-1}n\eta\max_{1\leq j\leq l}|S_{j}|_{1}
\right)
$$
avec $\eta=\frac{2l|L^{-1}|}{\min_k t_0^{\kappa_k}}$.

En conclusion, il existe $\eta_1>0$ tel que si $1/2<r\leq 1$, $m=2^k$, 
$\max_i\left|NF^{m+d_i-1}-NF^{p_0+d_i-1}\right|_{t_0^{\kappa}.r}<\eta_1$ et $\max_i\left|D\left(NF^{m+d_i-1}-J^1(NF^{m+d_i-1})\right)\right|_{t_0^{\kappa}.r}< \eta_1$, alors 
\begin{eqnarray}
|U_{\alpha}|_{t_0^{\kappa}.r} & \leq &\frac{c_1}{\omega_{k+1}(S)^2}\max_i|R_{i,\alpha}^{m+d_i}|_{t_0^{\kappa}.r}.
\end{eqnarray}

On peut supposer que $\frac{c_1}{\omega_{k+1}(S)^2}\geq 1$ quitte \`a prendre une plus grande valeur pour $c_1$, 
on peut donc \'ecrire $\frac{c_1}{\omega_{k+1}(S)^2}\=def \gamma_k^{-m}$
avec $\gamma_k\leq 1$. La majoration pr\'ec\'edente devient alors 
\begin{equation}\label{N-estim-gamma}
|U_{\alpha}|_{t_0^{\kappa}.r}\leq 
\gamma_k^{-m}\max_i|R_{i,\alpha}^{m+d_i}|_{t_0^{\kappa}.r}.
\end{equation}

\section{La r\'ecurrence}

Soient $1/2<r\leq 1$ un nombre r\'eel et $\eta_1>0$ le r\'eel positif d\'efinit par le th\'eor\`eme 
$(\ref{estim-cohom})$.
Pour tout entier $m\geq \max([8n/\min_k t_0^{\kappa_k}\eta_1]+1, p_0, \max d_i)$, on pose
\begin{eqnarray*}
{\cal NF}_{i,m}(r)& = & \left\{X\in  \pvf n {d_i} m\;|\; |X-J^{p_0+d_i-1}(X)|_{t_0^{\kappa}.r}<\eta_1-\frac{8n}{m-d_i+1},\right.\\
& & \left. |D(X-J^1(X))|_{t_0^{\kappa}.r}<\eta_1-\frac{8n}{\min_k t_0^{\kappa_k}(m-d_i+1)}\right\}\\
{\cal B}_{m+1}(r) & = & \left\{X\in \vfo n {m+1}\;|\; |X|_{t_0^{\kappa}.r}<1\right\}\\
\end{eqnarray*}
Si $m=2^k$ pour un entier $k\geq 1$, on d\'efinit
$$
\rho\=def m^{-1/m}r \quad\text{et}\quad R\=def\gamma_km^{-2/m}r
$$
o\`u le nombre $\gamma_k=\left(\frac{c_1}{\omega_{k+1}(S)^2}\right)^{-1/m}$ est d\'efinit par $(\ref{N-estim-gamma})$.

Il est clair que $m^{1/m}\geq 1$.
On peut supposer que les nombres $\omega_k$ sont plus petit que $1$; par cons\'equent, 
on a $\ln \omega_k < 1/m\ln m$ de sorte que $\ln \omega_k -2/m\ln m< -1/m\ln m<0$, c'est-\`a-dire $R<\rho<r\leq 1$.

Soit $\{X_1,\ldots X_l\}$ une alg\`ebre de type Cartan relativement \`a $S$, de champs de vecteur holomorphe. 
On suppose que $X_1$ est r\'egulier par rapport \`a $S$ et on suppose $X_i$ normalis\'e \`a l'ordre $m+d_i-1$ et on \'ecrit 
$$
X_i= NF^{m+d_i-1}_i + R_i^{m+d_i}
$$ 
o\`u $NF^{m+d_i-1}_i$ est une forme normale de degr\'e $m+d_i-1$ relativement \`a $s$, la partie lin\'eaire de $X_1$, et  
$R_i^{m+d_i}\in  \vfo n {m+d_i}$. Soit $U=\sum_{\alpha}U_{\alpha}$ la solution des \'equations (\ref{equ-cohom0}), 
la somme portant sur les poids non-nuls de $S$ dans $\pvf n {m+1} {2m}$.

On se propose de d\'emontrer le r\'esultat suivant~:
\begin{prop}\label{induction}
Avec les notations ci-dessus, supposons que, pour tout $i=1, \ldots, l$, on ait 
$( NF^{m+d_i-1}_i,R_{i,m+d_i})\in {\cal NF}_{i,m+d_i-1}(r)\times {\cal B}_{m+d_i}(r)$. Si $m$ est 
suffisamment grand (disons $m>m_0$ ind\'ependant de $r$), alors
\begin{enumerate}
\item $\Phi\=def (Id+U)^{-1}\in \text{Diff}_1(\Bbb C^n,0)$ est un diffeomorphisme tel que $D_{t_0^{\kappa}.R}\subset\Phi(D_{t_0^{\kappa}.\rho})$,
\item $\Phi^*X_i=NF^{2m+d_i-1}_i+R_i^{2m+d_i}$ est normalis\'e \`a l'ordre $2m+d_i-1$,
\item $(NF^{2m+d_i-1}_i,R_i^{2m+d_i})\in {\cal NF}_{i,2m+d_i-1}(R)\times {\cal B}_{2m+d_i}(R)$.
\end{enumerate}
\end{prop}

\subsection{Le diff\'eomorphisme normalisant}

Posons, pour $i=1,\ldots,l$,
$$
X_i= NF^{m+d_i-1}_i + B_i + C_i
$$ 
o\`u $NF^{m+d_i-1}_i\in \pvf n {d_i} {m+d_i-1}$, 
$B_i\in \pvf n {m+d_i} {2m+d_i-1}$ et $C\in \vfo n {2m+d_i}$.
On \'ecrit $B_i$ selon la d\'ecomposition de Fitting de ${\pvf n {m+d_i} {2m+d_i-1}}$ relativement \`a $S$ : 
$$
B_i= B_{i,0}+B^+_i\;\text{ avec }
B_{0,i}\in \pvf {n,0} {m+d_i} {2m+d_i-1}(S)\;\text{ et }\;
B^+_i\in \left(\pvf {n,\alpha} {m+d_i} {2m+d_i-1}\right)^+(S).
$$

Posons $\Phi^{-1}=\text{Id}+ U\in \text{Diff}_1(\Bbb C^n,0)$ ainsi que $x=\Phi^{-1}(y)$, on obtient 
\begin{equation}
D(\Phi^{-1})(y)\Phi^*(X_i)(y) =  X_i(\Phi^{-1}(y)),\quad i=1,\ldots, l\label{conj-equ}.
\end{equation}

Puisque $\Phi^*X_i$ est normalis\'e \`a l'ordre $2m+d_i-1$, \'ecrivons $\Phi^*X_i(y)=\text{NF}^{m+d_i-1}_i(y)+B'_i(y)+C'_i(y)$, avec 
$B'_i\in \pvf n {m+d_i} {2m+d_i-1}$ et $C'\in \vfo n {2m+d_i}$.
L'\'equation de conjugaison $(\ref{conj-equ})$ s'\'ecrit alors~:
\begin{eqnarray}
(\text{Id}+D(U)(y))(\text{NF}^{m+d_i-1}_i+B'_i+C'_i)(y) & = & (\text{NF}^{m+d_i-1}_i+B_i+C_i)(\Phi^{-1}(y))\label{eq-conj}\\
& = & \text{NF}^m(y)+ D(\text{NF}^m)(y)U(y) + B(y) \nonumber\\ 
&  & + \left( B_i(\Phi^{-1}(y))-B_i(y)\right) + C_i(\Phi^{-1}(y))\nonumber\\
&  & + \text{NF}^{m+d_i-1}_i(\Phi^{-1}(y))\nonumber \\
&  & -\left(\text{NF}^{m+d_i-1}_i(y)+ D(\text{NF}^{m+d_i-1}_i)(y)U(y)\right);\nonumber
\end{eqnarray}
ce que l'on peut r\'e\'ecrire 
\begin{eqnarray*}
C'_i(y)+\left(B'_i(y)-B_i(y)+[\text{NF}^{m+d_i-1}_i,U](y)\right) & = & \left( B_i(\Phi^{-1}(y))-B_i(y)\right) + C_i(\Phi^{-1}(y))\\
& & + \text{NF}^{m+d_i-1}_i(\Phi^{-1}(y)) \\
& & -\left(\text{NF}^{m+d_i-1}_i(y)+ D(\text{NF}^{m+d_i-1})(y)U(y)\right)\\
& & - D(U)(y)(B'_i+C'_i)(y)\\
& \=def &  D_i(y).\label{equ-jet}
\end{eqnarray*}

Puisque l'ordre de $\left( B_i(\Phi^{-1}(y))-B_i(y)\right)$ est sup\'erieur ou \'egal \`a l'ordre de 
$D(B_i)(y)U(y)$ et l'ordre de 
$\text{NF}^{m+d_i-1}_i(\Phi^{-1}(y))-\left(\text{NF}^{m+d_i-1}_i(y)+ D(\text{NF}^{m+d_i-1}_i)(y)U(y)\right)$ 
est sup\'erieur ou \'egal \`a l'ordre de $D^2(\text{NF}^{m+d_i-1}_i)(y)(U(y),U(y))$, alors 
l'ordre de $D_i(y)$ est sup\'erieur ou \'egal \`a $2m+d_i$. 
Par cons\'equent, on a $J^{2m+d_i-1}(B'_i(y)-B_i(y)+[\text{NF}^{m+d_i-1}_i,U])=0$; soit $B'_i=B_{i,0}$

\subsection{Calcul du reste}

Dans le but de majorer la norme de $C'_i$, on \'ecrit l'\'equation $(\ref{eq-conj})$ sous la forme suivante~:
\begin{eqnarray*}
C'_i(y) & = & \left( \text{NF}^{m+d_i-1}_i(\Phi^{-1}(y))-\text{NF}^{m+d_i-1}_i(y)\right) + (B_i+C_i)(\Phi^{-1}(y))\\
& & -B'_i(y)- D(U)(y)(\text{NF}^{m+d_i-1}_i+B'_i+C'_i)(y).
\end{eqnarray*}

Puisque $\text{NF}^{m+d_i-1}_i(\Phi^{-1}(y))-\text{NF}^{m+d_i-1}_i(y)=\int_0^1{D(\text{NF}^{m+d_i-1}_i)(y+tU(y))U(y)dt}$, 
nous utiliserons l'\'equation suivante 
\begin{eqnarray}
C'_i(y) & = & \int_0^1{D(\text{NF}^{m+d_i-1}_i)(y+tU(y))U(y)dt} + (B_i+C_i)(\Phi^{-1}(y))\label{reste}\\
& & -B'_i(y)- D(U)(y)(\text{NF}^{m+d_i-1}_i+B'_i+C'_i)(y).\nonumber
\end{eqnarray}

\subsection{Estimation du diff\'eomorphisme normalisant}

Soit $\Phi=(Id+U)^{-1}$ le diff\'eomorphisme normalisant. Par hypoth\`ese, 
$NF^{m+d_i-1}_i\in {\cal NF}_{i,m}(r)$; on peut donc appliquer le th\'eor\`eme $(\ref{estim-cohom})$~: 
$$
|U|_{t_0^{\kappa}.r}\leq \frac{c_1}{\omega_{k+1}(S)^2}\max_i |B^+_i|_{t_0^{\kappa}.r}.
$$
Puisque $B^+_i\prec \bar B^+_i + \bar B_{i,0}\prec \bar R_i^{m+d_i}$, on a  $|B^+|_{t_0^{\kappa}.r}<1$; 
il s'en suit que $|U|_{t_0^{\kappa}.r}\leq \gamma_k^{-m}$. 
\begin{lemms}\label{diffeo-rayon}
Sous les hypoth\`eses pr\'ec\'edentes et $m$ est suffisamment grand (disons $m>m_0$), alors 
pour tout $0<\theta\leq 1$ et tout entier $1\leq i\leq n$, on a 
$|y_i+\theta U_i(y)|_{t_0^{\kappa}.R}<t_0^{\kappa_i}\rho$ et
par cons\'equent, $\Phi(D_{t_0^{\kappa}.\rho})\supset D_{t_0^{\kappa}.R}$.
\end{lemms}
\Preuve 

Nous adaptons l'argument de Bruno \cite{Bruno}[p. 203].
Il suffit de montrer que $t_0^{\kappa_i}R+|U|_{t_0^{\kappa}.R}<t_0^{\kappa_i}\rho$.
Puisque $U$ est d'ordre $\geq m+1$ alors, par $(\ref{norm-rayon})$ et l'inegalit\'e 
pr\'ec\'edente, on a
\begin{eqnarray}
|U|_{t_0^{\kappa}.R} & \leq & \left(\frac{R}{r}\right)^{m+1}|U|_{t_0^{\kappa}.r}\nonumber\\
& \leq  & \left( \gamma_k m^{-2/m} \right)^{m+1}\gamma_k^{-m}\nonumber\\
& \leq & \gamma_k m^{-2-2/m} \label{u_k}\\
& \leq &  m^{-2-2/m}\nonumber.
\end{eqnarray}
Or $R = \gamma_k m^{-2/m}r\leq m^{-2/m}r$, il suffit donc de montrer que 
$m^{-2/m}(t_0^{\kappa_i}r+m^{-2})<t_0^{\kappa_i}\rho=t_0^{\kappa_i}m^{-1/m}r$; c'est-\`a-dire 
$\frac{m^{-2}}{m^{1/m}-1}<t_0^{\kappa_i}r$. Or,  
$$
\frac{m^{-2}}{m^{1/m}-1}=\frac{m^{-2}}{\exp^{1/m\ln m}-1}\leq \frac{m^{-2}}{1/m\ln m}\leq 
\frac{1}{m\ln m}
$$
car $1+x\leq\exp x$ pour tout $x\in \Bbb R^+$. Mais pour $0<x$ suffisamment grand, 
on a $2/\min_k t_0^{\kappa_k} < x\ln x$; par cons\'equent et puisque $1/2<r$, on obtient 
$\frac{m^{-2}}{m^{1/m}-1}< \frac{t_0^{\kappa_i}}{2}< t_0^{\kappa_i}r$.
\cqfd

\subsection{Majoration du reste}

On a $\left(\frac{R}{r}\right)^{m+1}=\gamma_k^{m+1}m^{-2-2/m}$, $\gamma_k\leq 1$ ainsi que 
$\frac{\rho}{r}=m^{-1/m}$; d'apr\`es ce qui pr\'ec\`ede, on
\begin{eqnarray}
|U|_{t_0^{\kappa}.R} & \leq & \gamma_km^{-2-2/m}\quad\text{par}\;(\ref{u_k}),\nonumber\\
& \leq & \frac{R}{r}m^{-2}\label{maj-U1}\\
& < &  m^{-2}\label{maj-U2}.
\end{eqnarray}

Comme nous l'avons vu, $\Phi^*(X_i)$ est normalis\'e \`a l'ordre $2m+d_i-1$ et 
$J^{2m+d_i-1}(\Phi^*(X_i))\=def NF^{2m+d_i-1}_i=NF^{m+d_i-1}_i+B_{i,0}$. 
Puisque $B_{i,0}$ est un polyn\^ome, il est domin\'e par $\bar B_i$, lequel est domin\'e par $\bar R_i^{m+d_i}$;
on a donc $|B_{i,0}|_{t_0^{\kappa}.r}\leq |R_i^{m+d_i}|_{t_0^{\kappa}.r}<1$.
On a alors
\begin{eqnarray*}
|B_{i,0}|_{t_0^{\kappa}.R} & \leq & (\gamma_k m^{-2/m})^{m+d_i}|B_{i,0}|_{t_0^{\kappa}.r},\\
& \leq & (\gamma_k m^{-2/m})^{m+1}|B_{i,0}|_{t_0^{\kappa}.r},\\
& \leq & m^{-2};\\
|D(B_{i,0})|_{t_0^{\kappa}.R} & \leq & \frac{2m+d_i-1}{\min_k t_0^{\kappa_k}R}|B_{i,0}|_{t_0^{\kappa}.R}\quad\text{par}\;(\ref{norm-diff}),\\
& \leq & \frac{2m+d_i-1}{\min_k t_0^{\kappa_k}R}\left(\frac{R}{r}\right)^{m+d_i}|B_{i,0}|_{t_0^{\kappa}.r}\\
& \leq & \frac{2m+d_i-1}{\min_k t_0^{\kappa_k}r}\left(\frac{R}{r}\right)^{m+d_i-1}|B_{i,0}|_{t_0^{\kappa}.r}\\
& \leq & \frac{2(2m+d_i-1)}{\min_k t_0^{\kappa_k}m^2}.
\end{eqnarray*}
On en d\'eduit que  
\begin{eqnarray*}
|NF^{2m+d_i-1}_i-NF^{p_0+d_i-1}_i|_{t_0^{\kappa}.R} & = & |(NF^{m+d_i-1}_i-NF^{p_0+d_i-1}_i)+B_{i,0}|_{t_0^{\kappa}.R} \\
& \leq & |NF^{m+d_i-1}_i-NF^{p_0+d_i}_i|_{t_0^{\kappa}.R}+|B_{i,0}|_{t_0^{\kappa}.R},\\
& < & \eta_1-\frac{8n}{m}+\frac{1}{m^2},\\
& < & \eta_1 -\frac{8n}{2m}\quad\text{si}\quad1<4nm;\\
|D(NF^{m+d_i-1}_i-J^{1}(NF^{m+d_i-1}_i))|_{t_0^{\kappa}.R} & = & |D(NF^{m+d_i-1}_i-J^{1}(NF^{m+d_i-1}_i))+D(B_{i,0})|_{t_0^{\kappa}.R} \\
& \leq & |D(NF^{m+d_i-1}_i-J^{1}(NF^{m+d_i-1}_i))|_{t_0^{\kappa}.R}+|D(B_{i,0})|_{t_0^{\kappa}.R},\\
& < & \eta_1-\frac{8n}{\min_k t_0^{\kappa_k}m} + \frac{2(2m+d_i-1)}{\min_k t_0^{\kappa_k}m^2} \\
& \leq & \eta_1-\frac{8n}{2m\min_k t_0^{\kappa_k}}\;\text{ si }\;(2n-2)m\geq (d_i-1).
\end{eqnarray*}
Ainsi $NF^{2m+d_i-1}_i\in {\cal NF}_{i,2m+d_i-1}(R)$. Il nous reste donc \`a montrer que $R_i^{2m+d_i}\in {\cal B}_{2m+d_i}(R)$.

On a les majorations suivantes~:
\begin{eqnarray*} 
|(B_i+C_i)\circ\Phi^{-1}|_{t_0^{\kappa}.R} & \leq & |B_i+C_i|_{t_0^{\kappa}.\rho} \quad\text{par}\;(\ref{diffeo-rayon}),\\
& \leq & \left(m^{-1/m}\right)^{m+d_i}|B_i+C_i|_{t_0^{\kappa}.r}\quad\text{car $B_i+C_i$ est d'ordre}\; \geq m+d_i,\\
& \leq & m^{-1} \\
|D(U)(NF^{m+d_i-1}_i+B_{i,0})|_{t_0^{\kappa}.R} & \leq & n|D(U)|_{t_0^{\kappa}.R}(|NF^{m+d_i-1}_i|_{t_0^{\kappa}.R}+|B_{i,0}|_{t_0^{\kappa}.R}),\\
& \leq & \frac{2nm}{\min_k t_0^{\kappa_k}R}|U|_{t_0^{\kappa}.R}(|NF^{m+d_i-1}_i|_{t_0^{\kappa}.R}+|B_{i,0}|_{t_0^{\kappa}.R})\\
& & \text{ car }\;U\;\text{est un polyn\^ome de degr\'e 2m},\\
& \leq & \frac{2nm}{\min_k t_0^{\kappa_k}r}m^{-2}(|NF^{m+d_i-1}_i|_{t_0^{\kappa}.R}+|B_{i,0}|_{t_0^{\kappa}.R}) 
\quad \text{par}\quad (\ref{maj-U1}),\\
& \leq & \frac{4n}{\min_k t_0^{\kappa_k}m}(|NF^{m+d_i-1}_i|_{t_0^{\kappa}.r}+m^{-2}) \quad\text{car}\quad r\geq 1/2;\\
|D(U)(y)C'_i|_{t_0^{\kappa}.R} & \leq & \frac{4n}{\min_k t_0^{\kappa_k}m}|C'_i|_{t_0^{\kappa}.R}\quad\text{par le m\^eme argument}.\\
\end{eqnarray*}

De plus, pour tout $0\leq \theta\leq 1$, on a
$$
D(NF^{m+d_i-1}_i)(y+\theta U(y))U(y)\prec D(\overline{NF^{m+d_i-1}_i})(y+ \bar U(y))\bar U(y);
$$
il s'en suit que  
\begin{eqnarray*}
\left|\int_0^1{D(NF^{m+d_i-1}_i)(y+\theta U(y))U(y)d\theta}\right |_{t_0^{\kappa}.R} & \leq & 
|D(\overline{NF^{m+d_i-1}_i})(y+ \bar U(y))\bar U(y)|_{t_0^{\kappa}.R},\\
&\leq & n|D(NF^{m+d_i-1}_i)|_{t_0^{\kappa}.\rho}|U|_{t_0^{\kappa}.R}\\
& \leq & \frac{n(m+d_i-1)}{\min_k t_0^{\kappa_k}\rho}|NF^{m+d_i-1}_i|_{t_0^{\kappa}.\rho}|U|_{t_0^{\kappa}.R}\\
& & \text{car }\;NF^{m+d_i-1}\;\text{est un polyn\^ome de degr\'e }m+d_i-1,\\
& \leq & \frac{n(m+d_i-1)}{\min_k t_0^{\kappa_k}r}|NF^{m+d_i-1}_i|_{t_0^{\kappa}.r}|U|_{t_0^{\kappa}.R}\\
& & \text{since}\;NF^{m+d_i-1}(0)=0,\\
& \leq & \frac{n(m+d_i-1)}{\min_k t_0^{\kappa_k}r}|NF^{m+d_i-1}_i|_{t_0^{\kappa}.r}\frac{R}{r}m^{-2}\\
& & \text{par }\;(\ref{maj-U1})\\
& \leq & 4n\frac{m+d_i-1}{\min_k t_0^{\kappa_k}m}|NF^{m+d_i-1}_i|_{r}m^{-1}\\
\end{eqnarray*}

On a donc, d'apr\`es $(\ref{reste})$
\begin{eqnarray*}
|C'_i|_{t_0^{\kappa}.R} & \leq & \left|\int_0^1{D(\text{NF}^{m+d_i-1}_i)(y+tU(y))U(y)dt}\right|_{t_0^{\kappa}.R} + |(B_i+C_i)(\Phi^{-1}(y))|_{t_0^{\kappa}.R}\\
& & +|B'_i|_{t_0^{\kappa}.R}+ |D(U)(y)(\text{NF}^{m+d_i-1}_i+B_{i,0})|_{t_0^{\kappa}.R}+|D(U)(y)C'_i|_{t_0^{\kappa}.R}\\
&\leq & 4n\frac{m+d_i-1}{\min_k t_0^{\kappa_k}m}|NF^{m+d_i-1}_i|_{t_0^{\kappa}.r}m^{-1}+m^{-1}+\frac{4n}{\min_k t_0^{\kappa_k}m}(|NF^{m+d_i-1}_i|_{t_0^{\kappa}.r}+m^{-2})\\
& & +\frac{4n}{\min_k t_0^{\kappa_k}m}|C'_i|_{t_0^{\kappa}.r}\\
& \leq & (\frac{m+d_i-1}{m}+)\frac{4n}{\min_k t_0^{\kappa_k}m}|NF^{m+d_i-1}_i|_{t_0^{\kappa}.r} + \frac{4n}{\min_k t_0^{\kappa_k}m^3} 
+ \frac{4n}{\min_k t_0^{\kappa_k}m}|C'_i|_{R}.
\end{eqnarray*}
Si $m>\frac{4n}{\min_k t_0^{\kappa_k}}$ alors 
\begin{eqnarray*}
|C'_i|_{t_0^{\kappa}.R} & \leq & \frac{4n}{\min_k t_0^{\kappa_k}m-4n}\left((1+\frac{m+d_i-1}{m})|NF^{m+d_i-1}_i|_{t_0^{\kappa}.r}+\frac{1}{m^2}\right)\\
& \leq & \frac{4n}{m-4n}((1+\frac{m+d_i-1}{m})(|NF^{p_0+d_i}_i|_{1}+\eta_1)+2).
\end{eqnarray*}
Soit $M$ la borne sup\'erieure des $\frac{m+d_i-1}{m}$; si $m>\frac{4n}{\min_k t_0^{\kappa_k}}\left((M(|NF^{p_0+d_i}_i|_1+\eta_1)+2)+1\right)$ alors $|C'_i|_{t_0^{\kappa}R}<1$, 
c'est-\`a-dire $C'_i=R_{i,2m+d_i}\in {\cal B}_i^{2m+d_i}(R)$.

\cqfd

\section{Preuve du th\'eor\`eme principal}

Nous reprenons l'argument classique qui permet de conclure, par r\'ecurrence. 
Soient $1/2< r\leq 1$ et la suite $\{R_k\}_{k\in \Bbb N}$ de r\'eels positifs d\'efinie par $R_0=r$, 
$R_{k+1}=\gamma_km^{-2/m}R_k$ avec $m=2^k$. On a le 
\begin{lemm}\cite{Bruno, Stolo-sintg2}\label{nb}
La suite $\{R_k\}_{k\in \Bbb N}$ converge et il existe un entier $k_1$ tel que $\forall k\geq k_1$, $R_k>R_{k_1}/2$.
\end{lemm}

Soientt alors $X_1,\ldots, X_l$ comme dans l'\'enonc\'e. On peux supposer que l'alg\`ebre de type Cartan est normalis\'ee 
\`a un ordre $m_2=2^{k_2}$ suffisamment grand ($\geq \max(m_0, 2^{k_1}$) et on \'ecrit 
$X_i=NF^{m_2+d_i-1}_i+R_i^{m_2+d_i}$, $i=1,\ldots, l$.
Quitte \`a faire agir une homoth\'etie, on peut supposer que 
$(NF^{m_2+d_i-1}_i,R_i^{m_2+d_i})\in {\cal NF}_{i,m_2}(1)\times {\cal B}_{m_2+d_i}(1)$. On red\'efinit la suite 
$\{R_k\}$ par $R_{k_2}=1$. En vertu du lemme \ref{nb}, on a pour tout $k\geq k_0$, $R_k>1/2$.

Montrons par r\'ecurrence sur $k\geq k_2$, qu'il existe un diff\'eomorphisme $\Psi_k$ de $(\Bbb C^n,0)$ tel que 
pour tout $i=1,\ldots,l$, $\Psi_k^*(NF^{m_2+d_i-1}_i+R^{m_2+d_i}_i)\=def NF^{2^{k_2+1}+d_i-1}_i+R^{2^{k+1}+d_i}_i$ soit 
normalis\'e \`a l'ordre $2^{k+1}+d_i-1$, 
$(NF^{2^{k+1}+d_i-1}_i,R^{2^{k+1}+d_i}_i)\in {\cal NF}_{i,2^{k+1}+d_i-1}(R_{k+1})\times {\cal B}_{2^{k+1}+d_i}(R_{k+1})$ et
$|\text{Id} -\Psi_k^{-1}|_{t_0^{\kappa}.R_{k+1}}\leq \sum_{p=k_0}^k\frac{1}{2^{2p}}$.

\begin{itemize}
\item Pour $k=k_2$: 
D'apr\`es la proposition $(\ref{induction})$, il existe un diff\'eomorphisme $\Phi_{k_2}$ tel que 
$\Phi_{k_2}^*(NF^{m_2+d_i-1}_i+R^{m_2+d_i}_i)=NF^{2m_2+d_i-1}_i+R^{2m_2+d_i}_i$ est normalis\'e \`a 
l'ordre $2m_2+d_i-1$, $(NF^{2m_2+d_i-1}_i,R^{m_2+d_i})\in {\cal NF}_{i,2m_2+d_i-1}(R_{k_2+1})\times {\cal B}_{2m_2+d_i}(R_{k_2+1})$ et 
$|\text{Id} -\Phi_{k_2}^{-1}|_{t_0^{\kappa}.R_{k_2+1}}<1/2^{2{k_2}}$.

\item Supposons que le r\'esultat soit vrai pour tout entier $i\leq k-1$ : par hypoth\`eses, 
$\Psi_{k-1}^*(NF^{m_2+d_i-1}_i+R^{m_2+d_i}_i)= NF^{2^{k}+d_i-1}_i+R^{2^{k}+d_i}_i$ est normalis\'e \`a l'ordre $2^{k}$ et
$(NF^{2^{k}+d_i-1}_i,R^{2^{k}+d_i}_i)\in {\cal NF}_{i,2^{k}+d_i-1}(R_{k})\times {\cal B}_{2^{k}+d_i}(R_{k})$.
Puisque $1/2<R_k\leq 1$, on peut appliquer la proposition $(\ref{induction})$ : il existe un diff\'eomorphisme $\Phi_{k}$ 
tel que $(\Phi_{k}\circ\Psi_{k-1})^*(NF^{m_2+d_i-1}_i+R^{m_2+d_i}_i)=NF^{2^{k+1}+d_i-1}_i+R^{2^{k+1}+d_i}_i$ est
normalis\'e \`a l'ordre $2^{k+1}+d_i-1$ et 
$(NF^{2^{k+1}+d_i-1}_i,R^{2^{k+1}+d_i}_i)\in {\cal NF}_{i,2^{k+1}+d_i-1}(R_{k+1})\times {\cal B}_{2^{k+1}+d_i}(R_{k+1})$.
Posons $\Psi_{k}=\Phi_{k}\circ \Psi_{k-1}$. D'apr\`es la proposition $(\ref{induction})$ (ou le lemme $(\ref{diffeo-rayon})$), 
on a $|\text{Id} - \Phi_{k}^{-1}|_{t_0^{\kappa}.R_{k+1}}<1/2^{2{k}}$. Il s'en suit que 
\begin{eqnarray*}
|\text{Id} -\Psi_{k}^{-1}|_{t_0^{\kappa}.R_{k+1}}& \leq & \left|(\text{Id}-\Psi_{k-1}^{-1})\circ \Phi_{k}^{-1}+(\text{Id} -\Phi_{k}^{-1})\right|_{t_0^{\kappa}.R_{k+1}},\\
& \leq & \left|(\text{Id}-\Psi_{k-1}^{-1})\circ \Phi_{k}^{-1}\right|_{t_0^{\kappa}.R_{k+1}}+\left|(\text{Id} -\Phi_{k}^{-1})\right|_{t_0^{\kappa}.R_{k+1}}.\\
\end{eqnarray*}
D'apr\`es la proposition $(\ref{induction})$, on a $\Phi_{k}^{-1}(D_{t_0^{\kappa}.R_{k+1}})\subset D_{t_0^{\kappa}.R_k}$; 
par cons\'equent, 
\begin{eqnarray*}
|\text{Id} -\Psi_{k}^{-1}|_{t_0^{\kappa}.R_{k+1}}& \leq & \left|(\text{Id}-\Psi_{k-1}^{-1})\right|_{t_0^{\kappa}.R_{k}}+\left|(\text{Id} -\Phi_{k}^{-1})\right|_{t_0^{\kappa}.R_{k+1}}\\
& \leq & \sum_{p=k_0}^{k-1}\frac{1}{2^{2p}}+\frac{1}{2^{2k}};
\end{eqnarray*}
ce qui termine la preuve de la r\'ecurrence.
\end{itemize}
Puisque $D_{t_0^{\kappa}.1/2}\subset D_{t_0^{\kappa}.R_{k}}$ pour tout entier $k\geq k_2$, la suite $\{|\Psi_k^{-1}|_{t_0^{\kappa}.1/2}\}_{k\geq k_2}$ est 
uniform\'ement born\'ee. De plus, la suite $\{\Psi_k^{-1}\}_{k\geq k_2}$ converge vers le diff\'eomorphisme 
formel $\hat\Psi^{-1}$ (l'inverse du diff\'eomorphisme normalisant) dans l'espace des s\'eries formelles. 
Par cons\'equent, cette suite converge 
dans ${\cal H}_n^n(t_0^{\kappa}.r)$ (pour tout $r<1/2$) vers $\hat\Psi^{-1}$ (c.f. \cite{Grauert-L1}). Cela signifie 
que la transformation normalisante est holomorphe au voisinage de $0\in \Bbb C^n$.

\section{Applications}
\subsection{Le cas hamiltonien d'Ito}

Soit $H_1=\sum_{k=1}^n{\lambda_k x_ky_k}+\tilde H$ une fonction holomorphe au voisinage de $0\in \Bbb C^{2n}$; 
$\lambda_i\in \Bbb C$, $\tilde H$ \'etant d'ordre $\geq 3$ en $0$. On suppose que 
$$
\forall (m_1,\ldots, m_n)\in \Bbb Z^n\setminus\{0\},\quad \sum_{k=1}^n\lambda_km_k\neq 0.\quad\quad (*)
$$
On consid\`ere le syst\`eme d'\'equations diff\'erentielles 
\begin{equation}
\frac{dx_k}{dt}=\frac{\partial H_1}{\partial y_k},\quad \frac{dy_k}{dt}=-\frac{\partial H_1}{\partial x_k},\quad k=1,\ldots,n \label{hamilt}
\end{equation}
H. Ito a d\'emontr\'e les r\'esultat suivant, qui g\'en\'eralise le travail de J. Vey dans le cadre hamiltonien \cite{vey-ham}~:
\begin{theo}\cite{Ito1}
Sous les hypoth\`eses pr\'ec\'edentes, on suppose que le syst\`eme (\ref{hamilt}) admet $n-1$ autres int\'egrales premi\`eres 
$H_2,\ldots, H_n$ holomorphes au voisinage de $0\in \Bbb C^{2n}$;  et on suppose que $H_1,\ldots, H_n$ sont 
fonctionnellement ind\'ependantes. Il existe alors une transformation $\phi$ canonique (i.e hamiltonienne) et
holomorphe au voisinage de $0\in \Bbb C^{2n}$ telle que $H_k\circ \phi$ soient une fonction holomorphe des $n$ mon\^omes 
$x_ky_k$.
\end{theo}
\Preuve
Soient $\lie g$ une alg\`ebre de Lie complexe commutative de dimension $n$, $\{g_1,\ldots,g_n\}$ une base de $\lie g$ et 
$S$ le morphisme de Lie d\'efinit par $S(g_i)=x_i\frac{\partial}{\partial x_i}-y_i\frac{\partial}{\partial y_i}$. 
Le morphisme $S$ est diophantien car les poids de $S$ dans les espaces de champs homog\`enes prennent de valeurs 
enti\`eres sur la base de $\lie g$.
On pose $X_i$ les champs de vecteurs hamiltoniens d\'efinit par les fonctions $H_i$. Gr\^ace \`a la condition 
$(*)$ la partie lin\'eaire de $X_1$ est r\'eguli\`ere par rapport \`a $S$. Le fait qu'il soit hamiltonien implique 
que $X_1$ est r\'egulier par rapport \`a $S$. D'autres par, les champs $X_i$, $i\geq 2$, commutent avec $X_1$ et sont hamiltoniens.
Par cons\'equent, dans un bon syst\`eme $\Phi$ de coordonn\'ees formelles commun, ils appartiennent au $\widehat{\cal O}_n^S$-module 
engendr\'e par $S(\lie g)$. Gr\^ace au lemme de Ziglin, on peut choisir eventuellement d'autres int\'egrales premi\`eres de sorte 
que leur ind\'ependance fonctionelle implique que 
la libert\'e des parties juniors de leur champs associ\'es sur $\widehat{\cal O}_n^S$. On remarquera que $\widehat{\cal O}_n^S=\Bbb C[[x_1y_1,\ldots, x_ny_n]]$. D'apr\`es le 
th\'eor\`eme \ref{th-princ}, $\Phi$ est holomorphe au voisinage de $0\in \Bbb C^{2n}$. Un argument du \`a J. Vey, permet 
alors de modifier $\Phi$ en un diff\'eomorphisme symplectique holomorphe tout en changant de forme normale.
\cqfd

\bibliographystyle{alpha}
\bibliography{normal,math,asympt,analyse,stolo,lie}

\begin{thebibliography}{PMY94}

\bibitem[Arn80]{Arn2}
V.~Arnold.
\newblock {\em {Chapitres suppl\'{e}mentaires de la th\'{e}orie des
  \'{e}quations diff\'{e}rentielles ordinaires}}.
\newblock Mir, 1980.

\bibitem[Bru72]{Bruno}
A.D Bruno.
\newblock {The analytical form of differential equations}.
\newblock {\em Trans. Mosc. Math. Soc}, 25,131-288(1971); 26,199-239(1972),
  1971-1972.

\bibitem[CS]{Camacho-port}
C.~Camacho and P.~Sad.
\newblock {Pontos singulares de equa\c{c}oes diferencais analiticas}.
\newblock 16 Coloquio Brasileiro de Matematica.

\bibitem[Eca92]{Ecalle-mart}
J.~Ecalle.
\newblock {Singularit\'es non abordables par la g\'eom\'etrie}.
\newblock {\em Ann. Inst. Fourier, Grenoble}, 42,1-2(1992),73-164, 1992.

\bibitem[GR71]{Grauert-L1}
H.~Grauert and R.~Remmert.
\newblock {\em Analytische Stellenalgebren}.
\newblock Springer-Verlag, 1971.

\bibitem[Ito89]{Ito1}
H.~Ito.
\newblock Convergence of birkhoff normal forms for integrable systems.
\newblock {\em Comment. Math. Helv.}, 64:412--461, 1989.

\bibitem[Mal82]{Malg-diffeo}
B.~Malgrange.
\newblock {Travaux d'Ecalle et de Martinet-Ramis sur les syst\`{e}mes
  dynamiques}.
\newblock {\em S\'{e}m. Bourbaki 1981-1982}, exp. 582(1982),, 1982.

\bibitem[MR82]{Ram-Mart1}
J.~Martinet and J.P. Ramis.
\newblock {Probl\`{e}mes de modules pour des \'{e}quations diff\'{e}rentielles
  non lin\'{e}aires du premier ordre}.
\newblock {\em I.H.E.S}, 55,63-164, 1982.

\bibitem[MR83]{Ram-Mart2}
J.~Martinet and J.P. Ramis.
\newblock {Classification analytique des \'{e}quations diff\'{e}rentielles non
  lin\'{e}aires r\'{e}sonantes du premier ordre}.
\newblock {\em Ann. Sci. E.N.S}, 4\`{e}me s\'{e}rie,16,571-621, 1983.

\bibitem[PMY94]{perez-yoccoz}
R.~P\'erez-Marco and J.-C. Yoccoz.
\newblock Germes de feulletages holomorphes \`a holonomie pr\'escrite.
\newblock In {\em Complex analytic methods in dynamical systems}, volume 222 of
  {\em Ast\'erisque}, pages p.345--371. Soc. Math. France, 1994.

\bibitem[Rou75]{roussarie-ast}
R.~Roussarie.
\newblock Mod\`{e}les locaux de champs et de formes.
\newblock {\em Ast\'{e}risque}, 30, 1975.

\bibitem[Sie42]{Siegel}
C.L. Siegel.
\newblock {Iterations of analytic functions}.
\newblock {\em Ann. Math.}, 43(1942)807-812, 1942.

\bibitem[Sto96]{Stolo-clas}
L.~Stolovitch.
\newblock {Classification analytique de champs de vecteurs $1$-r\'{e}sonnants
  de $(\Bbb C^n,0)$}.
\newblock {\em Asymptotic Analysis}, 12:91--143, 1996.

\bibitem[Sto98]{Stolo-intg-cras}
L.~Stolovitch.
\newblock Compl\`{e}te int\'{e}grabilit\'{e} singuli\`{e}re.
\newblock {\em C.R. Acad. Sci., Paris, S\'{e}rie I}, 326:733--736, 1998.

\bibitem[Sto99]{Stolo-sintg2}
L.~Stolovitch.
\newblock Singular complete integrabilty.
\newblock Technical Report 142, Pr\'{e}publication E. Picard, Janvier 1999.
\newblock 1-67.

\bibitem[Sto00]{Stolo-cartan-cras}
L.~Stolovitch.
\newblock Normalisation holomorphe d'alg\`ebres de type cartan de champs de
  vecteurs holomorphes singulier.
\newblock {\em C.R. Acad. Sci, Paris, S\'{e}rie I}, 330:1--4, 2000.

\bibitem[Vey78]{vey-ham}
J.~Vey.
\newblock Sur certains syst\`{e}mes dynamiques s\'{e}parables.
\newblock {\em Am. Journal of Math. 100}, pages 591--614, 1978.

\bibitem[Vor81]{Voronin}
S.M. Voronin.
\newblock {Analytic classification of germs of conformal mappings $(\Bbb
  C,0)\rightarrow (\Bbb C,0)$ with identity linear part}.
\newblock {\em Funct. An. and its Appl.}, 15,(1981), 1981.

\bibitem[Yoc88]{Yoccoz}
J.-C. Yoccoz.
\newblock {Lin\'{e}arisation des germes de diff\'{e}omorphismes holomorphes de
  $(\Bbb C, 0)$}.
\newblock {\em C.R. Acad. Sci. Paris}, 306(1988)55-58, 1988.

\bibitem[Yoc95]{Yoccoz-ast}
J.-C. Yoccoz.
\newblock Petits diviseurs en dimension 1.
\newblock {\em Ast\'{e}risque}, 231, 1995.

\end{thebibliography}
\end{document}